\documentclass[12pt]{iopart}
\usepackage{iopams}
\usepackage{graphicx}
\usepackage{latexsym}

\begin{document}
\newtheorem{theorem}{Theorem}[section]
\newtheorem{lemma}[theorem]{Lemma}
\newtheorem{corollary}[theorem]{Corollary}
\newtheorem{Remark}[theorem]{Remark}
\newenvironment{remark}{\begin{Remark}\rm}{\end{Remark}}
\newtheorem{Definition}[theorem]{Definition}
\newenvironment{definition}{\begin{Definition}\rm}{\end{Definition}}
\newenvironment{proof}
    {\rm \trivlist \item[\hskip \labelsep{\bf Proof. }]}
    {\hspace*{\fill}$\Box$\endtrivlist}

\eqnobysec

\def\diag{\textrm{diag}}
\def\ds{\displaystyle}

\title[The relativistic Toda lattice and Laurent orthogonal polynomials]
{Direct and inverse spectral transform for the relativistic Toda
lattice and the connection with Laurent orthogonal polynomials}

\author{J Coussement, A B J Kuijlaars and W Van Assche}

\address{Department of Mathematics, Katholieke Universiteit Leuven,
Celestijnenlaan 200B, 3001 Leuven, Belgium}


\begin{abstract}
We introduce a spectral transform for the finite relativistice
Toda lattice (RTL) in generalized form.  In the nonrelativistic
case, Moser constructed a spectral transform from the spectral
theory of symmetric Jacobi matrices.  Here we use a non-symmetric
generalized eigenvalue problem for a pair of bidiagonal matrices
$(L,M)$ to define the spectral transform for the RTL.  The inverse
spectral transform is described in terms of a terminating
T-fraction.  The generalized eigenvalues are constants of motion
and the auxiliary spectral data have explicit time evolution.
Using the connection with the theory of Laurent orthogonal
polynomials, we study the long-time behaviour of the RTL.  As in
the case of the Toda lattice the matrix entries have asymptotic
limits.  We show that $L$ tends to an upper Hessenberg matrix with
the generalized eigenvalues sorted on the diagonal, while $M$
tends to the identity matrix.
\end{abstract}

\section{Introduction}

The relativistic Toda lattice (RTL) was introduced by
Ruijsenaars \cite{Rui} and studied in \cite{Brus1,Brus2,Cosentino,Ragnisco,Suris},
for a review see \cite{Kharchev}. The finite RTL is defined by the
system of equations
\begin{equation}
\fl \label{Newtoniaans} \qquad \qquad \ddot q_{n} =
\epsilon^2\dot q_{n} \left(\dot
q_{n-1}\frac{\exp(q_{n-1}-q_{n})}{1+\epsilon^2\exp(q_{n-1}-q_{n})}-\dot
q_{n+1}\frac{\exp(q_{n}-q_{n+1})}{1+\epsilon^2\exp(q_{n}-q_{n+1})}\right),
\end{equation}
with $N\in \mathbb{N}$, $1 \leq n \leq N$ and the convention $q_{0} \equiv -\infty$
and $q_{N+1} \equiv +\infty$.  Here the $q_n$ are functions in the
time parameter $t$.  The system (\ref{Newtoniaans}) is the
Newtonian form of a Hamiltonian system with Hamiltonian
\[ H(q_{1},\ldots,q_{N},p_{1},\ldots,p_{N}) =
    \sum_{n=1}^N\rme^{p_{n}}h(q_{n-1}-q_{n})h(q_{n}-q_{n+1}), \]
where $h(x)=\sqrt{1+\epsilon^2\rme^x}$.  Setting $\dot q_{n}=\dot
Q_{n}+1/\epsilon$ and letting $\epsilon\to 0$, one can easily
check that  (\ref{Newtoniaans}) reduces to the equations of motion
of the finitie nonrelativistic Toda lattice,
\[ \ddot Q_{n}=\exp(Q_{n-1}-Q_{n})-\exp(Q_{n}-Q_{n+1}), \qquad 1\le n\le N, \]
so that (\ref{Newtoniaans}) is a one-parameter perturbation.

In analogy with the Flaschka variables for the nonrelativistic
Toda lattice one can use a (non-invertible) change of variables to
prove integrability. Bruschi and Ragnisco \cite{Brus1,Brus2}, see
also \cite{Kharchev,Suris}, obtain the two forms
\begin{equation}
\label{hetsysteem} \left\{ \begin{array}{ll} \dot
a_{n}=\frac{b_{n}}{a_{n+1}}-\frac{b_{n-1}}{a_{n-1}}, & \qquad 1\le n\le N,
\\[2ex] \dot
b_{n}=b_{n}\left(\frac{1}{a_{n}}-\frac{1}{a_{n+1}}\right),& \qquad
1\le n\le N-1,
\end{array}\right.
\end{equation}
 and
\begin{equation}
\label{anderevorm} \left\{ \begin{array}{ll} \dot
a_{n}= a_{n}(b_{n-1}-b_{n}), & \qquad 1\le n\le N,\\[2ex]
\dot b_{n}= b_{n}(a_{n}-a_{n+1}+a_{n-1}-a_{n+1}),&
\qquad 1\le n\le N-1,
\end{array}\right.
\end{equation}
both with $a_n > 0$ for $1 \leq n \leq N$, $b_n > 0$ for $1 \leq n
\leq N-1$, and $b_0 \equiv 0$, $b_{N}\equiv 0$. The systems
(\ref{hetsysteem}) and (\ref{anderevorm}) can also be written in
matrix form. Define two bidiagonal matrices $L$ and $M$ by
\begin{eqnarray} \nonumber L=\left(
\begin{array}{cccccc}
a_{1} & 1 &  0 & \cdots & \cdots & 0\\ 0 & a_{2} & 1 &  & & \vdots
\\ 0 & 0 & a_3 &  \ddots &  &  \vdots\\ \vdots & & \ddots & \ddots
& \ddots & \vdots \\ \vdots & & & \ddots & a_{N-1}& 1 \\ 0 &
\cdots & \cdots & \cdots & 0 & a_{N}
\end{array} \right), \\
\label{dematricesLenU}
M=\left(
\begin{array}{cccccc}
1 & 0 & 0 & \cdots & \cdots & 0\\ -b_{1} & 1 & 0 &  & & \vdots \\
0 & -b_2 & 1 &  \ddots &  &  \vdots\\ \vdots & & \ddots & \ddots &
\ddots & \vdots\\ \vdots & & & \ddots & 1& 0 \\ 0 & \cdots &
\cdots &\cdots& -b_{N-1} & 1
\end{array} \right).
\end{eqnarray}
Then Suris \cite{Suris} noted that (\ref{hetsysteem}) can be
written in the Lax form
\begin{equation} \label{Laxform}
    \left\{\begin{array}{l}
    \dot L=LA-BL\\[2ex]
    \dot M=MA-BM,
    \end{array}\right.
\end{equation}
 where
\begin{equation} \fl \label{matricesAenB}
\qquad \qquad A=\left(
\begin{array}{ccccc}
0 & \cdots & \cdots & \cdots & 0\\ \frac{b_{1}}{a_{2}} & 0 & & &
\vdots \\ 0 & \ddots & \ddots &   &  \vdots\\ \vdots & \ddots &
\ddots & \ddots & \vdots\\ 0 & \cdots & 0& \frac{b_{N-1}}{a_{N}} &
0
\end{array} \right),\quad B=\left(
\begin{array}{ccccc}
0 & \cdots & \cdots & \cdots & 0\\ \frac{b_{1}}{a_{1}} & 0 & & &
\vdots \\ 0 & \ddots & \ddots &   &  \vdots\\ \vdots & \ddots &
\ddots & \ddots & \vdots\\ 0 & \cdots & 0& \frac{b_{N-1}}{a_{N-1}}
& 0
\end{array} \right).
\end{equation}
Then $A=-(L^{-1}M)_{-}$ and $B=-(ML^{-1})_{-}$, where we use
$X_{-}$ to denote the strictly lower triangular part of $X$. The
other system (\ref{anderevorm}) can also be written in the form
(\ref{Laxform}) but now $A = -(M^{-1}L)_-$ and $B= -(LM^{-1})_-$,
see Remark \ref{stellinganderevorm} below. In Section 3.1 we show
that (\ref{hetsysteem}) and (\ref{anderevorm}) are special cases
of a generalized form of the finite RTL, namely
\begin{equation}
\label{moregeneralintro} \left\{\begin{array}{l} \dot
L=\Bigl(F\left(LM^{-1}\right)\Bigr)_{-}\
L-L\Bigl(F\left(M^{-1}L\right)\Bigr)_{-}\\ \dot
M=\Bigl(F\left(LM^{-1}\right)\Bigr)_{-}\
M-M\Bigl(F\left(M^{-1}L\right)\Bigr)_{-}
\end{array}\right.
\end{equation}
where $F : \mathbb R^+ \to \mathbb R$ is an arbitrary function on $\mathbb{R}^{+}$.
For $F(x) = 1/x$, we obtain (\ref{hetsysteem}) and
for $F(x) = x$, we obtain (\ref{anderevorm}).

\begin{figure}[h]
\begin{center}
\setlength{\unitlength}{1truemm}
\begin{picture}(120,70)(0,0)
\put(0,0){\framebox(30,15){\parbox{2.5cm}{\begin{center}spectral
data\\ at $t=0$\end{center}}}} \put(32,7.5){\vector(1,0){56}}
\put(90,0){\framebox(30,15){\parbox{2.5cm}{\begin{center}spectral
data\\ at $t$\end{center}}}}
\put(0,55){\framebox(30,15){\parbox{2.5cm}{\begin{center}matrix
data\\ at $t=0$\end{center}}}}
\put(90,55){\framebox(30,15){\parbox{2.5cm}{\begin{center}matrix
data\\ at $t$\end{center}}}} \put(32,62.5){\vector(1,0){56}}
\put(15,53){\vector(0,-1){36}} \put(105,17){\vector(0,1){36}}
\put(17,35){\parbox{2cm}{\begin{center}direct spectral
problem\end{center}}}
\put(83,35){\parbox{2cm}{\begin{center}inverse spectral
problem\end{center}}}
\put(40,13.5){\parbox{4cm}{\begin{center}evolution of the spectral
data\end{center}}}
\end{picture}
\end{center}
\caption{\label{schemaplot}{\em The scheme of the direct and
inverse spectral problem.}}
\end{figure}
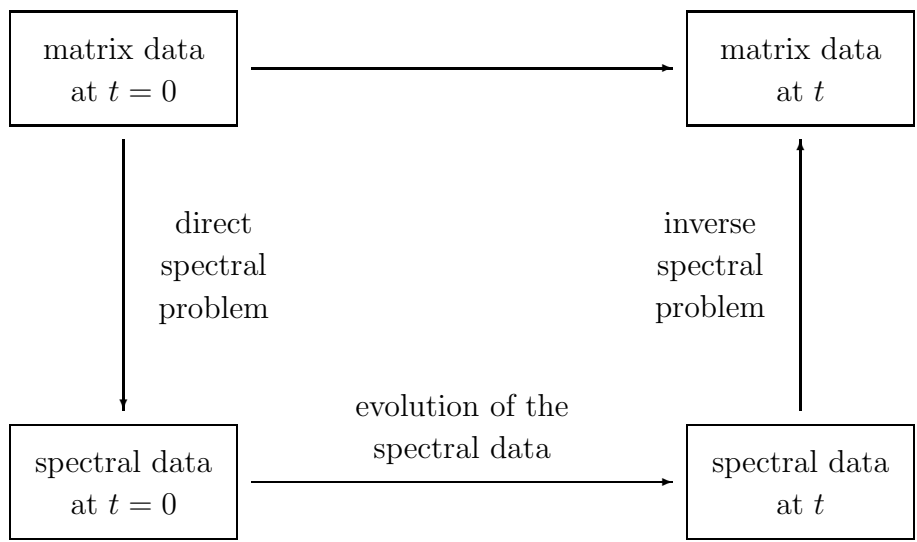

The main goal of this paper is to solve the finite RTL in its
general form (\ref{moregeneralintro}) with the aid of a direct and
inverse spectral transform, in a similar way as was done by Moser
\cite{Moser} for the nonrelativistic case. While the spectral
transform for the nonrelativistic case is based on the spectral
theory of tridiagonal (Jacobi) matrices and is connected with
orthogonal polynomials, we will show that in the relativistic
case, the spectral transform is based on the spectral theory of
pairs of bidiagonal matrices (\ref{dematricesLenU}) and Laurent
orthogonal polynomials \cite{Jones4,Ranga,Vinet,Hendrik}. See
Figure \ref{schemaplot} for the general scheme of the direct and
inverse spectral problem.

To the best of our knowledge, the only transform methods that have
been explicitly introduced in the study of the RTL are the ones by
Ragnisco and Bruschi \cite{Ragnisco} for the periodic RTL,  by
Ruijsenaars \cite{Rui}, and the scattering transform of Cosentino
\cite{Cosentino} for the infinite RTL.

The Laurent orthogonal polynomials $P_n$ are defined by
\begin{equation} \label{Laurentpol}
    P_0\equiv 1,\qquad P_n(z)=\det(zM_n-L_n), \qquad  1\le n\le N,
\end{equation}
where $L_n$ and $M_n$ are the $n\times n$
upper left corner blocks of $L$ and $M$, respectively.
They satisfy the recurrence relation
\[ P_{n}(z) = (z-a_{n}) P_{n-1}(z) - b_{n-1}z P_{n-2}(z), \qquad 1\le n\le N, \]
with $P_0\equiv 1$, $P_{-1}\equiv 0$.  Under the assumption that
$a_n
> 0$ and $b_n
> 0$, it is known that the zeros of these polynomials are simple,
real and positive, see \cite{Jones3,Jones4,Jones,Ranga} and Lemma
\ref{interlacing} below. The first set of spectral data we choose
are the ordered zeros of the polynomial $P_N(z)=\det(zM-L)$ which
we denote by $0<\lambda_1<\cdots<\lambda_N$. Then the $\lambda_j$
are also the eigenvalues for the generalized eigenvalue problem of
the pair $(L,M)$
\begin{equation} \label{eigproblem}
    L \vec{x} = \lambda M \vec{x}.
\end{equation}
To obtain our second set of spectral data we
associate with each $\lambda_j$ the value
\[w_j= \left( \vec{u}_j^{\, T} M \vec{v}_j \right)^{-1}, \]
where $\vec{u}_j$ and $\vec{v}_j$ are the left and right
eigenvectors of the generalized eigenvalue problem
(\ref{eigproblem}), respectively, normalized to have their first
component equal to 1. We prove in Section \ref{dsp} that $w_j > 0$
and $\sum_{j=1}^N w_j = 1$. The direct spectral transform
\[ \left\{ \begin{array}{cc}a_1, \ldots, a_N  > 0 \\
           b_1, \ldots, b_{N-1} > 0 \end{array} \right\} \longmapsto
    \left\{ \begin{array}{cc}
    0 < \lambda_1 < \cdots < \lambda_N \\
    w_1, \ldots, w_N > 0, \ \sum_{j=1}^N w_j = 1 \end{array} \right\} \]
is invertible.  Using $T$-fractions \cite{Hendrik,Jones}, we
explain in Section \ref{isp} how to obtain the $a_n$ and $b_n$
when the spectral data $\lambda_j$ and $w_j$ are given. We also
note that the polynomials $P_n$ satisfy
\[ \sum_{j=1}^N P_n(\lambda_j) P_m(\lambda_j) \frac{w_j}{\lambda_j^n} = 0,
    \qquad \mbox{if } n > m, \]
which is known as Laurent orthogonality
\cite{Cochran,Jones3,Jones4,Ranga}.

In Section \ref{section_evolution} we describe the time evolution
of the spectral data for the finite RTL in the generalized form
(\ref{moregeneralintro}) with positive initial data. We show that
the $\lambda_j$ are time independent and that
\[
w_{j}(t)=\frac{w_{j}(0)\rme^{-F\left(\lambda_{j}\right)t}}{\sum\limits_{k=1}^Nw_{k}(0)\rme^{-
F\left(\lambda_{k}\right)t}},\qquad  1\le j\le N.
\]

Finally, in Section \ref{section_limits} we investigate the
behaviour of the RTL for $t \to \pm\infty$. We assume that $F$ is
strictly increasing. Then
\begin{equation} \label{limitbn}
    \lim_{t \to \pm \infty} b_n(t) = 0
\end{equation}
and
\begin{equation} \label{limitan}
    \lim_{t\to -\infty} a_{n}(t)= \lambda_{N-n+1},
    \qquad \lim_{t\to +\infty}a_{n}(t)=\lambda_{n},
    \qquad 1\le n\le N.
\end{equation}
Similar limits hold for the case that $F$ is strictly decreasing.
This sorting property of the RTL indicates that the choice of the
spectral data $\lambda_1 < \cdots < \lambda_N$ is very natural. In
Theorems \ref{blimieten} and \ref{alimieten} we establish the
precise rates of convergence of (\ref{limitbn}) and
(\ref{limitan}). This is similar to Deift et al. \cite{Deift} who
obtained precise rates of convergence for the SVD flow (Singular
Value Decomposition).

\section{The direct and inverse spectral problem}

\label{detransformatie} In this section we introduce a
transformation from matrix data $a_n>0$, $b_n>0$ to spectral data
$\lambda_j$, $w_j$. We prove that this spectral transform  is a
bijection between $\left(\mathbb{R}^+\right)^{2N-1}$ and its image
$\{ 0 < \lambda_1 < \cdots < \lambda_N, \ w_j >0, \ \sum_{j=1}^N
w_j = 1\}$.

\subsection{The direct spectral problem}
\label{dsp}
We start from positive matrix data $a_n$, $1 \leq n \leq N$,
and $b_{n}$, $1 \leq n \leq N-1$, and construct the bidiagonal matrices
$L$ and $M$ as in (\ref{dematricesLenU}).

Define the finite set of monic polynomials $P_n$ as
\[ P_0\equiv 1,\qquad P_{n}(z)=\det(zM_n-L_n), \qquad 1\le n\le  N \]
where $L_n$ and $M_n$ are obtained from the matrices $L$ and $M$,
respectively, by deleting the last $N-n$ rows and columns.
These polynomials satisfy the recurrence relation
\begin{equation}
\label{recursie1}
P_{n}(z)=(z-a_{n})P_{n-1}(z)-b_{n-1}zP_{n-2}(z),\qquad 1\le n\le N,
\end{equation}
as can be easily seen by expanding the determinant
$\det(zM_n-L_n)$ first by the last column and then by the last
row. The initial conditions for (\ref{recursie1}) are $P_{0}
\equiv 1$ and $P_{-1} \equiv 0$. The zeros of these polynomials
behave in a way which is similar to the  behaviour of zeros of
polynomials which are orthogonal on the positive real line. The
following lemma can be found in \cite{Jones3,Jones4,Jones,Ranga},
but we include a proof for completeness.
\begin{lemma}
\label{interlacing} The zeros of the polynomials
$P_n(z)=\det(zM_n-L_n)$ are real, simple, positive
($0<z_{1,n}<\cdots<z_{n,n}$) and have the interlacing property,
which means that
\[ 0 < z_{1,n}<z_{1,n-1}<z_{2,n}<z_{2,n-1}<\cdots<z_{n-1,n-1}<z_{n,n}. \]
\end{lemma}
\begin{proof}  The lemma is true for the case $n=1$ because $P_1(z)=z-a_1$
and $a_1>0$. Suppose that the lemma is true for some $n < N$. By
the interlacing property we then know that $P_{n-1}$ changes sign
between any two consecutive zeros of $P_{n}$. Evaluating the
recurrence
\[ P_{n+1}(z)=(z-a_{n+1})P_{n}(z)-b_{n}zP_{n-1}(z)
\]
at a zero $z_{j,n}$ of $P_n$, we get $P_{n+1}(z_{j,n}) =
-b_{n}z_{j,n}P_{n-1}(z_{j,n})$. Since $b_{n}>0$ and $z_{j,n}>0$,
we see that $P_{n+1}$ and $P_{n-1}$ have opposite signs at the
zeros of $P_n$. It follows that also $P_{n+1}$ changes sign
between any two consecutive zeros of $P_{n}$, and therefore it has
at least one zero in each of the intervals $(z_{j,n}, z_{j+1,n})$,
$1\le j\le n-1$. Since $P_{n-1}(z_{n,n})$ is positive, it also
follows that $P_{n+1}(z_{n,n})$ is negative. Since $P_{n+1}$ is a
monic polynomial, it must have a zero in $(z_{n,n}, \infty)$.

To complete the proof we will show that there is also a zero in
$(0,z_{1,n})$. Evaluating the recurrence at $0$, we get
$P_{n+1}(0) = - a_{n+1} P_{n}(0)$. Since $P_n$ is a monic
polynomial of degree $n$ having all its zeros on the positive real
line, we find that $(-1)^n P_{n}(0) > 0$. Thus $(-1)^n
P_{n+1}(0) < 0$. Since $P_{n-1}$ is a monic polynomial of degree
$n-1$ with all its zeros to the right of $z_{1,n}$, we obtain
$(-1)^n P_{n-1}(z_{1,n}) < 0$. Since $P_{n+1}$ and $P_{n-1}$
differ in sign at $z_{1,n}$, we get $(-1)^n P_{n+1}(z_{1,n}) > 0$.
So $P_{n+1}$ has opposite signs in $0$ and $z_{1,n}$, which shows
that there is indeed a zero in $(0,z_{1,n})$. The lemma now
follows by mathematical induction.
\end{proof}

\begin{remark}
\label{deltapolynomen} An analogous reasoning as in Lemma
\ref{interlacing} can be applied to the polynomials
$\Delta_n(z)=\det (z\tilde M_n-\tilde L_n)$, $\Delta_{0}\equiv 1$,
where $\tilde L_n$ and $\tilde M_n$ are obtained from the matrices
$L$ and $M$ by deleting the {\em first} $N-n$ rows and columns,
respectively. These polynomials satisfy the recurrence
\begin{equation}
\label{recursie2}
\Delta_{n}(z)=(z-a_{N-n+1})\Delta_{n-1}(z)-b_{N-n+1}z\Delta_{n-2}(z),
\end{equation}
for $1 \leq n \leq N$,
with initial conditions $\Delta_{0} \equiv 1$ and $\Delta_{-1} \equiv 0$.
So also the zeros of these polynomials are real, simple and
positive and have the interlacing property.
\end{remark}

Consider now the generalized eigenvalue problem of the pair
$(L,M)$
\begin{equation}
\label{gep} L\vec{x}=\lambda M\vec{x}.
\end{equation}
The eigenvalues are the zeros of the polynomial
$P_N(\lambda)=\det(\lambda M-L)$.  From Lemma \ref{interlacing} we
know that these eigenvalues are real, simple and positive. We
denote them by
\[ 0<\lambda_{1}<\cdots<\lambda_{N}<\infty. \]
We use $\vec{v}_j$ to denote the right eigenvector corresponding
to the eigenvalue $\lambda_j$, normalized so that the first
component is equal to $1$. From the bidiagonal structure of the
matrices $M$ and $L$ it easily follows that the first component of
any eigenvector of (\ref{gep}) is non-zero, so that this
normalization is always possible. We use $\vec{u}_j$ to denote the
corresponding left eigenvector, also normalized to have the first
component equal to $1$. We view $\vec{u}_j$ as a column vector so
that $\vec{u}_j^{\, T}L=\lambda_j \vec{u}_j^{\, T}M$. We collect
the eigenvectors in the matrices
\[V = \left(\vec{v}_1 \cdots  \vec{v}_N \right)\qquad \mbox{and}
    \qquad U = \left(\vec{u}_1 \cdots  \vec{u}_N \right). \]
A left and right eigenvector of this generalized eigenvalue
problem, corresponding to different eigenvalues, are
$M$-orthogonal, that is, $\vec{u}_j^{\, T}M\vec{v}_k=0$ whenever
$j \neq k$. This implies that $U^TMV$ is a diagonal matrix. Since
the three matrices $U$, $M$, and $V$ are invertible, the diagonal
elements of $U^TMV$ are non-zero, and we can make the following
definition.

\begin{definition}
\label{definitieW_j}
 Let $U$ and $V$ be as above. Then we define
\begin{equation}
\label{expressionw_j} w_j=\frac{1}{(U^TMV)_{j,j}} = \Bigl(
\vec{u}_j^{\, T}M\vec{v}_j \Bigr)^{-1},\qquad 1\le j\le  N.
\end{equation}
\end{definition}

Now we have defined the spectral data $\lambda_j$ and $w_j$ for
$1\le j\le N$. To complete the description of the direct spectral
transform we show that $w_j>0$ and $\sum_{j=1}^Nw_j=1$. To prove
this we introduce the Weyl (or resolvent) function
\[ f(z)=\vec{e}_1^{\, T}(zM-L)^{-1}\vec{e}_1=\Bigl((zM-L)^{-1}\Bigr)_{1,1},
\qquad z\in \mathbb{C}\setminus \{\lambda_1,\ldots,\lambda_N\},\]
with $\vec{e}_1$ the first unit vector in $\mathbb{R}^N$.
\begin{lemma}
\label{Wsd} The Weyl function is given by
\[ f(z)=\sum_{j=1}^N\frac{w_{j}}{z-\lambda_{j}}. \]
\end{lemma}
\begin{proof}
Note that by definition of $w_j$, we have $\vec{u}_j^{\, T} M
\vec{v}_k = \frac{1}{w_j} \delta_{j,k}$, and so
\[\vec{u}^{\, T}_j L \vec{v}_k = \lambda_j \vec{u}^{\, T}_j M\vec{v}_k
    = \frac{\lambda_j}{w_j}\delta_{j,k},
    \qquad 1\le j,k\le N. \]
In matrix form this is $U^{\,
T}MV=\diag\Bigl(\frac{1}{w_{1}},\ldots,\frac{1}{w_{N}}\Bigr)$ and
$U^{\,
T}LV=\diag\Bigl(\frac{\lambda_1}{w_{1}},\ldots,\frac{\lambda_N}{w_{N}}\Bigr).$
Combining these two equations, we find that
\[ U^T (zM-L) V = \diag\Bigl(\frac{z-\lambda_1}{w_{1}},\ldots,\frac{z-\lambda_N}{w_{N}}\Bigr).\]
From this we readily obtain that
\begin{equation} \label{Weyl2}
    (zM-L)^{-1}= V \diag\Bigl(\frac{w_{1}}{z-\lambda_1},\ldots,\frac{w_{N}}{z-\lambda_N}\Bigr)
    U^T.
\end{equation}
We now recall that the eigenvectors $\vec{u}_j$ and $\vec{v}_j$
are normalized with their first component equal to $1$. Thus the
first rows of the matrices $U$ and $V$ contain only ones. Then
(\ref{Weyl2}) gives that
\[f(z)=\Bigl((zM-L)^{-1}\Bigr)_{1,1}=\sum_{j=1}^N\frac{w_{j}}{z-\lambda_{j}}.\]
\end{proof}

The lemma shows that the $\lambda_j$ are simple poles of the Weyl
function and that the $w_j$ are the corresponding residues.  A
second representation of $f$ is in terms of the polynomials
$\Delta_n$ introduced in Remark \ref{deltapolynomen}.
\begin{lemma}
\label{WpT0} The Weyl function $f$ can be written as
$f(z)= \ds \frac{\Delta_{N-1}(z)}{\Delta_N(z)}$.
\end{lemma}
\begin{proof} This follows immediately from Cramer's rule.
\end{proof}

As a corollary we get
\begin{corollary}
We have $w_j > 0$ for every $j$ and $\sum_{j=1}^N w_j = 1$.
\end{corollary}
\begin{proof}If we combine Lemma \ref{Wsd} and Lemma
\ref{WpT0}, we get
\begin{equation}
\label{connectie}
\frac{\Delta_{N-1}(z)}{\Delta_N(z)}=\sum_{j=1}^N\frac{w_{j}}{z-\lambda_{j}}.
\end{equation}
From Remark \ref{deltapolynomen} we know that the zeros of
$\Delta_{N-1}$ and $\Delta_N$ are interlacing. Since
$\Delta_{N-1}$ and $\Delta_N$ are monic polynomials, it then
follows that $w_j > 0$ for every $j$. The fact that $\sum_{j=1}^N
w_j=1$  follows by multiplying (\ref{connectie}) by $z$ and then
letting $z \to \infty$.
\end{proof}

We now have the direct spectral transform
\[ \left( \mathbb R^+ \right)^{2N-1} \to \Lambda :
  (a_{1},\ldots,a_{N},b_{1},\ldots,b_{N-1})  \mapsto  (\lambda_{1},\ldots,\lambda_{N},w_{1},\ldots,w_{N})
\]
where
\begin{equation} \label{Lambda}
    \Lambda =\left\{(\lambda_{1},\ldots,\lambda_{N},w_{1},\ldots,w_{N})
\left|\ \begin{array}{cc} 0<\lambda_{1}<\cdots<\lambda_{N}<\infty, \\[8pt]
    w_{j}>0,\ \sum_{j=1}^Nw_{j}=1\end{array}\right.  \right\}.
\end{equation}

We conclude this paragraph by establishing an explicit
representation of the Weyl function in terms of the given matrix
data.  This representation is also the key to the inverse spectral
transform. It turns out that the function $zf(z)$ can be written
as a special continued fraction, known as a T-fraction, see
\cite{Jones,Hendrik}.
\begin{lemma}
\label{WpT} The Weyl function $f$ can be written as
\begin{equation} \label{Tfraction}
f(z)=\frac{\Delta_{N-1}(z)}{\Delta_N(z)}=\frac{1}{z-a_{1}-
\frac{\displaystyle{b_{2}z}}{\displaystyle{z-a_{2}-}
\frac{\displaystyle{b_{3}z}}{
\begin{array}{ccc}
z-a_{3}- & & \\
 & \ddots & \\
 & & -\frac{\displaystyle{b_{N}z}}{\displaystyle{z-a_{N}}}
\end{array}
}}}.
\end{equation}
\end{lemma}
\begin{proof}
Let $r_k(z) = \frac{\Delta_{k-1}(z)}{\Delta_k(z)}$ for $1\le k\le
N$. Since $\Delta_0(z) = 1$ and $\Delta_1(z) = z-a_N$ we have
$r_1(z) = \frac{1}{z-a_{N}}$. From the recurrence
(\ref{recursie2}) satisfied by the polynomials $\Delta_k$, we get
that
\[ r_{k}(z)=\frac{1}{z-a_{N-k+1}-zb_{N-k+2}r_{k-1}(z)}. \]
Then it is easy to show by induction that
\begin{equation} \label{rkz} \fl
r_k(z) = \frac{\Delta_{k-1}(z)}{\Delta_k(z)}
    = \frac{1}{z-a_{N-k+1}-
\frac{\displaystyle{b_{N-k+2}z}}{\displaystyle{z-a_{N-k+2}-}\frac{\displaystyle{b_{N-k+3}z}}{
\begin{array}{ccc}
z-a_{N-k+3}- & & \\
 & \ddots & \\
 & & -\frac{\displaystyle{b_{N}z}}{\displaystyle{z-a_{N}}}
\end{array}
}}}.
\end{equation}
Taking $k=N$ in (\ref{rkz}) we obtain the lemma.
\end{proof}

From Lemma \ref{Wsd} and Lemma \ref{WpT} it is easy to see that
the spectral transform is injective.  In the next section we show
that it is also surjective.

\subsection{The inverse spectral problem}
\label{isp} We now start from spectral data
$(\lambda_{1},\ldots,\lambda_{N},w_{1},\ldots,w_{N})$ in
$\Lambda$, see (\ref{Lambda}), and show that there exist positive
matrix data $(a_{1},\ldots,a_{N},b_{1},\ldots,$ $b_{N-1})$ which
correspond to the given spectral data. We also give a method to
construct these matrix data.  In view of Lemma \ref{WpT} it is our
task to develop $\sum_{j=1}^N \frac{w_j}{z-\lambda_j}$ into a
continued fraction of the form (\ref{Tfraction}) with positive
$a_n$ and $b_n$, that is, we want
\begin{equation}
\label{claim}
\sum_{j=1}^N\frac{w_{j}}{z-\lambda_{j}}=\frac{1}{z-a_{1}-
\frac{\displaystyle{b_{1}z}}{\displaystyle{z-a_{2}-}
\frac{\displaystyle{b_{2}z}}{
\begin{array}{ccc}
z-a_{3}- & & \\
 & \ddots & \\
 & & -\frac{\displaystyle{b_{N-1}z}}{\displaystyle{z-a_{N}}}
\end{array}
}}}.
\end{equation}

\begin{theorem}
\label{existence_of_matrix_data} For given spectral data
$(\lambda_{1},\ldots,\lambda_{N},w_{1},\ldots,w_{N})$ in $\Lambda$
there exist positive $a_n$ and $b_n$ such that {\rm (\ref{claim})}
holds.
\end{theorem}
\begin{proof}Let us say that $f\in {\cal W}_N$ if $f$ is of the
form
\begin{equation}
\label{vormvanf} f(z)=\sum_{j=1}^N\frac{w_{j}}{z-\lambda_{j}},
\end{equation}
with $0 < \lambda_1 < \cdots < \lambda_N$, $w_j > 0$, and $\sum_{j=1}^Nw_{j}=1$.
We will prove that for
$f\in {\cal W}_N$, there exist $a, b> 0$ and $g\in {\cal W}_{N-1}$ so that
\[ f(z)=\frac{1}{z-a-bzg(z)}. \]
The theorem then follows by repeated application of this fact.

To prove the claim, we observe that $f$ is a rational function
which we can write as
\[f(z)=\frac{p_{N-1}(z)}{q_{N}(z)},\]
with $p_{N-1}$ and $q_{N}$ monic polynomials of degrees $N-1$ and $N$, respectively
(here we use the fact that $\sum_{j=1}^Nw_{j}=1$).  Then there exist $a, b \in \mathbb R$
and a monic polynomial $p_{N-2}$ of degree at most $N-2$ such that
\begin{equation}
\label{kiesaenb} q_N(z)=(z-a)p_{N-1}(z)-bzp_{N-2}(z).
\end{equation}
Taking $z=0$ in (\ref{kiesaenb}) we find
\begin{equation}
\label{formulevoora}
    a=-\frac{q_N(0)}{p_{N-1}(0)}=
    \left(\sum\limits_{j=1}^N\frac{w_{j}}{\lambda_{j}}\right)^{-1},
\end{equation}
which is positive, since all $\lambda_j$ and $w_j$ are positive.
To find $b$, we write (\ref{kiesaenb}) as
$-q_N(z) + (z-a)p_{N-1}(z) = bzp_{N-2}(z)$ and compute the coefficient
of $z^{N-1}$ on the left-hand side. Since $q_N(z) = \prod_{j=1}^N (z-\lambda_j)$
and $p_{N-1}(z) = \sum_{k=1}^N w_k \prod_{j \neq k} (z-\lambda_j)$, this
coefficient is
\begin{eqnarray*}
\sum_{j=1}^N\lambda_{j}-\sum_{k=1}^Nw_{k} \sum_{j=1,j\neq k}^N \lambda_{j} - a & = &
    \sum_{j=1}^N\lambda_{j}\left(1-\sum_{k=1,k\neq j}^Nw_{k}\right)-a\\
 & = & \sum_{j=1}^Nw_{j}\lambda_{j}-\left(\sum_{j=1}^N\frac{w_j}{\lambda_j}\right)^{-1}.
\end{eqnarray*}
In the last step we used (\ref{formulevoora}) and $\sum_{j=1}^N
w_j = 1$. Next, from Jensen's inequality
\[ \phi \left(\int z \rmd\mu(z)\right)< \int \phi(z) \rmd\mu(z) \]
applied to the strictly convex function $\phi(z) = 1/z$ and the
probability measure  $\mu = \sum_{j=1}^Nw_j\delta_{\lambda_j}$, we
get that
\begin{equation} \label{waardevanb}
\sum_{j=1}^Nw_{j}\lambda_{j}-\left(\sum_{j=1}^N\frac{w_j}{\lambda_j}\right)^{-1}
> 0.
\end{equation}
So the coefficient of $z^{N-1}$ in $bz p_{N-2}(z)$ is positive.
Since $p_{N-2}$ is a monic polynomial of degree $\leq N-2$, it
then follows that $p_{N-2}$ has exact  degree $N-2$ and that $b$
is equal to (\ref{waardevanb}). Hence $b$ is positive. If we now
define $g=p_{N-2}/p_{N-1}$ then we have that
\[f(z)=\frac{1}{z-a-bzg(z)},\qquad \mbox{ with } a, b > 0. \]
It remains to show that $g\in {\cal W}_{N-1}$.
Observe that the poles of $g$ are the zeros of $f$.
From (\ref{vormvanf}) with $w_j > 0$ it follows that
the zeros of $f$ interlace with its poles
$0 < \lambda_1 < \lambda_2 < \cdots < \lambda_N$.
So the poles of $g$ are simple and positive. If we denote
them by
$0 < \lambda_1^{\star} < \lambda_2^{\star} < \cdots < \lambda_{N-1}^{\star}$,
then
\[ g(z)=\frac{p_{N-2}(z)}{p_{N-1}(z)} =
    \sum_{j=1}^{N-1}\frac{w^{\star}_{j}}{z-\lambda^{\star}_{j}},\]
for certain residues $w^{\star}_j$ that satisfy
$\sum_{j=1}^{N-1}w^{\star}_{j}=1$, because $p_{N-2}$ and $p_{N-1}$ are monic.
We also have
\begin{eqnarray*}
w^{\star}_j &  = &
\lim \limits_{z \rightarrow \lambda^{\star}_{j}}(z-\lambda^{\star}_{j})g(z)
 \ = \ \lim \limits_{z \rightarrow \lambda^{\star}_{j}}\frac{z-
 \lambda^{\star}_{j}}{b}\left(\frac{z-a}{z}-\frac{1}{zf(z)}\right) \\
 & = & -\lim \limits_{z \rightarrow \lambda^{\star}_{j}}
 \frac{z-\lambda^{\star}_{j}}{bzf(z)}
 \ =  \ \frac{-1}{b\lambda^{\star}_{j}f'(\lambda^{\star}_{j})}.
\end{eqnarray*}
Since $b > 0$, $\lambda^{\star}_j > 0$, and $f'(\lambda^{\star}_j)
= - \displaystyle \sum_{k=1}^N
\frac{w_k}{(\lambda^{\star}_j-\lambda_k)^2} < 0$, we see that
$w^{\star}_j > 0$ for every $j$.\\ This completes the proof of the
theorem.
\end{proof}

\subsection{Laurent orthogonality}
Suppose that $\mu$ is a positive measure with support in $\mathbb
R^+$, for which all the strong moments $\int z^k \rmd\mu(z)$, $k
\in \mathbb{Z}$, exist. For $n$ not exceeding the number of mass
points in the support of $\mu$, the monic Laurent orthogonal
polynomial, $p_n$, of degree $n$ is defined as the monic
orthogonal polynomial of degree $n$ for the varying measure
$d\mu_n(z)=z^{-n}d\mu(z)$. So $p_n$ satisfies the orthogonality
relations
\[\int p_n(z)z^k \frac{\rmd\mu(z)}{z^n}=0,\qquad 0\le k\le n-1. \]
Cochran and Cooper \cite{Cochran} and Jones et al.
\cite{Jones2,Jones3,Jones4}  showed that these polynomials satisfy
a recurrence relation of the form
\[p_{n}(z)=(z-\alpha_n)p_{n-1}(z)-\beta_{n-1}zp_{n-2}(z),\qquad n\ge 1,\]
with initial conditions $p_0 \equiv 1$ and $p_{-1} \equiv 0$,
and positive recurrence coefficients $\alpha_n$ and $\beta_{n-1}$.
From the theory of orthogonal polynomials it follows that the zeros of
$p_n$ are simple, real and positive, since $\mu$ is a measure on the
positive real line.  Some examples of
(continuous) Laurent orthogonal polynomials are given in
\cite{Ranga_alleen}.

We now show that the polynomials $P_n(z)=\det (zM_n-L_n)$
introduced in (\ref{Laurentpol}) are the monic Laurent orthogonal polynomials
with respect to the discrete probability measure
\[\mu_N=\sum_{j=1}^Nw_{j}\delta_{\lambda_{j}}, \]
see also \cite{Ranga} for a similar situation. These polynomials
are also generated by the recurrence relation (\ref{recursie1}).
This illustrates the connection between the spectral transform,
proposed in Section \ref{dsp}, and the theory of Laurent
orthogonal polynomials.

\begin{theorem}
\label{Laurent-orthogonality}
 Let $a_{n}>0$, $1\le n\le N$, and $b_{n}>0$, $1\le n\le N-1$.
 Let the polynomials $P_n$ be defined by the recurrence relation {\rm (\ref{recursie1})}
 with $P_0 \equiv 1$ and $P_{-1} \equiv 0$. Let $\lambda_j$ and $w_j$
 be the spectral data associated with $a_n,b_n$ as defined in Section {\rm \ref{dsp}}. Then
\[ \sum_{j=1}^N
P_m(\lambda_j)P_n(\lambda_j)\frac{w_j}{\lambda_j^n}=\delta_{m,n}\prod_{k=1}^nb_k,\qquad
0\le m\le n\le N-1.
\]
\end{theorem}
\begin{proof} Evaluating the recurrence relation
(\ref{recursie1}) at $\lambda_j$ (which is a zero of $P_N$) one
can easily see that the right eigenvector $\vec{v}_j$ for the
generalized eigenvalue problem (\ref{gep}) is equal to
\begin{equation}
\label{expressionrighteigenvector}
    \vec{v}_j=\Bigl(1,P_1(\lambda_j),\ldots,P_{N-1}(\lambda_j)\Bigr)^T.
\end{equation}
Also note that the left eigenvector $\vec{u}_j$ can be written as
\begin{equation}
\label{expressionlefteigenvector}
\vec{u}_j=\left(1,\frac{P_1(\lambda_j)}{\lambda_jb_1},
\frac{P_2(\lambda_j)}{\lambda_j^2b_2b_1},\ldots,
\frac{P_{N-1}(\lambda_j)}{\lambda_j^{N-1}b_{N-1}\ldots
b_1}\right)^T.
\end{equation}
Since $U^{\,
T}MV=\diag\Bigl(\frac{1}{w_{1}},\ldots,\frac{1}{w_{N}}\Bigr)$ by
the definition of the $w_j$, we have
\[ V\diag(w_{1},\ldots,w_{N}) U^T = M^{-1}. \]
Now note that $M$ is lower triangular with ones on the diagonal,
and therefore the same holds for $M^{-1}$. So if $0\le m \leq n \leq N-1$
we have
\[\Bigl(V\diag(w_{1},\ldots,w_{N}) U^T\Bigr)_{m+1,n+1}=\delta_{m,n}.\]
Using (\ref{expressionrighteigenvector}) and
(\ref{expressionlefteigenvector}) we finally get
\[\sum_{j=1}^N
P_m(\lambda_j)P_n(\lambda_j)\frac{w_j}{\lambda_j^n \prod\limits_{k=1}^nb_k}
    =\delta_{m,n},\]
which proves the theorem.
\end{proof}

The theorem suggests an alternative method for performing the
inverse spectral transform. Knowing the spectral data, one can
construct the monic Laurent orthogonal polynomials corresponding
to the discrete measure $\mu_N$. These polynomials satisfy a three
term recurrence relation with coefficients equal to $a_n$ and
$b_n$.

\section{The generalized finite RTL and spectral evolution}

In the introduction we have already mentioned that the finite RTL
\begin{equation*}
\left\{ \begin{array}{ll} \dot
a_{n}=\frac{b_{n}}{a_{n+1}}-\frac{b_{n-1}}{a_{n-1}}, & \qquad 1\le
n\le N,
\\[2ex] \dot
b_{n}=b_{n}\left(\frac{1}{a_{n}}-\frac{1}{a_{n+1}}\right), &
\qquad 1\le n\le N-1,
\end{array}\right.
\end{equation*}
with $b_{0}\equiv 0$, $b_{N}\equiv 0$, can be written in  Lax form
representation \cite{Suris}, namely
\begin{equation}
\label{vormmetAenB} \left\{\begin{array}{l} \dot L=LA-BL\\[2ex]
\dot M=MA-BM,
\end{array}\right.
\end{equation}
with $A=-(L^{-1}M)_{-}$ and $B=-(ML^{-1})_{-}$. Here the matrices
$L$ and $M$ contain the matrix data $a_n$ and $b_n$ and are
defined in (\ref{dematricesLenU}).    In this section we first
show that this system can be generalized as in
(\ref{moregeneralintro}). Then we will solve the finite RTL in its
generalized form, using the spectral transform introduced in
Section \ref{detransformatie}. We will find the explicit evolution
of the spectral data.

\subsection{The generalized finite RTL}
\label{section_generalRTL}
The generalized finite RTL is defined
by the matrix differential equations
\begin{equation}
\label{moregeneral} \left\{\begin{array}{l} \dot
L=\Bigl(F\left(LM^{-1}\right)\Bigr)_{-}\
L-L\Bigl(F\left(M^{-1}L\right)\Bigr)_{-}\\ \dot
M=\Bigl(F\left(LM^{-1}\right)\Bigr)_{-}\
M-M\Bigl(F\left(M^{-1}L\right)\Bigr)_{-}\end{array}\right.
\end{equation}
where $F : \mathbb R^+ \to \mathbb R$ is a real-valued function on
$\mathbb R^+$.   The matrices $L$ and $M$ contain the matrix
data $a_n$ and $b_n$ as in (\ref{dematricesLenU}). We assume
positive initial data $a_n(0) > 0$, $b_n(0)>0$.  The system is of
the form (\ref{vormmetAenB}) if we set
\begin{equation}
\label{matrixAenBalgemeen}
A=-\Bigl(F\left(M^{-1}L\right)\Bigr)_{-}\qquad \mbox{and}\qquad
B=-\Bigl(F\left(LM^{-1}\right)\Bigr)_{-}.
\end{equation}

The eigenvalues of $LM^{-1}$ and $M^{-1}L$ are equal to the
generalized eigenvalues for the pair $(L,M)$. We know that there
are $N$ distinct positive generalized eigenvalues. Hence the
matrices $LM^{-1}$ and $M^{-1}L$ are diagonalizable, and in fact
we have
\[ LM^{-1} = U^{-T} D U^T, \qquad M^{-1}L = V^{-1} D V \]
where $D$ is a diagonal matrix containing the
generalized eigenvalues $\lambda_1 < \lambda_2 < \cdots < \lambda_N$.
Thus
\[ F\left(L M^{-1}\right) = U^{-T} \diag\Bigl(F(\lambda_{1}),\ldots,F(\lambda_{N}) \Bigr) U^T,
\]and\[ F\left(M^{-1}L\right) = V^{-1}
\diag\Bigl(F(\lambda_{1}),\ldots,F(\lambda_{N}) \Bigr)V,
\]
see \cite[Definition 6.2.4]{Horn}.

Of course one has to check that the structure of the matrices $L$ and
$M$ is preserved by the differential equations
(\ref{moregeneral}). We will use the following lemma.
\begin{lemma}
\label{UenLequatilties} For every function $F$ we have that
\begin{equation}
\label{Fisveelterm}
LF\left(M^{-1}L\right)=F\left(LM^{-1}\right)L\quad
\mbox{and} \quad
MF\left(M^{-1}L\right)=F\left(LM^{-1}\right)M.\end{equation}
\end{lemma}
\begin{proof}For a monomial $F(x) = x^n$, we have
\[LF\left(M^{-1}L\right) = L\left(M^{-1}L\right)^n =
\left(LM^{-1}\right)^n L = F\left(LM^{-1}\right)L\] and
\[MF\left(M^{-1}L\right)= M \left(M^{-1} L\right)^n =
    \left(LM^{-1}\right)^{n}M=F\left(LM^{-1}\right)M.\]
Then, by linearity, (\ref{Fisveelterm}) holds for every polynomial
$F$. On the finite set of eigenvalues of $LM^{-1}$ (which are the
same as the eigenvalues of $M^{-1}L$) each function $F$ is equal
to a polynomial, so the equalities hold for every $F$.
\end{proof}

We now prove that the proposed generalisation of the finite RTL is
justified.

\begin{theorem}
The system of matrix differential equations {\rm (\ref{moregeneral})} is
well defined, which means that $LA-BL$ is a diagonal matrix and
$MA-BM$ has only non-zero elements on the first subdiagonal, where
\[ A=-\Bigl(F\left(M^{-1}L\right)\Bigr)_{-} \quad \mbox{ and }
\quad
 B=-\Bigl(F\left(LM^{-1}\right)\Bigr)_{-}. \]
\end{theorem}
\begin{proof} Since $L$ is upper Hessenberg and $A$ and $B$ are
strictly lower triangular, one immediately sees that $LA-BL$ is
lower triangular.  Now $L$ is also upper triangular, so that
\[\left(\Bigl(F\left(LM^{-1}\right)\Bigr)_{-}\
L\right)_{-}=\Bigl(F\left(LM^{-1}\right)\ L\Bigr)_{-}\] and
\[\left(L\Bigl(F\left(M^{-1}L\right)\Bigr)_{-}\right)_{-}=
\Bigl(LF\left(M^{-1}L\right)\Bigr)_{-}.\] This means
that\begin{eqnarray*}
 (LA-BL)_{-}&= &
\Bigl(F\left(LM^{-1}\right)L-LF\left(M^{-1}L\right)\Bigr)_{-}.
\end{eqnarray*}
From Lemma \ref{UenLequatilties} we then obtain that
$(LA-BL)_{-}$ is the zero matrix.  We have
already mentioned that $LA-BL$ is lower triangular, so we conclude
that $LA-BL$ is a diagonal matrix.

Since $M$ is lower triangular
and $A$ and $B$ are strictly lower triangular, $MA-BM$
is strictly lower triangular.  Then we note that only the first of
the subdiagonals of $M$ has non-zero elements.  This implies that
below the first subdiagonal,
\[MA-BM=\Bigl(F\left(LM^{-1}\right)\Bigr)_{-}\
M-M\Bigl(F\left(M^{-1}L\right)\Bigr)_{-}\] agrees with
\[F\left(LM^{-1}\right)M-MF\left(M^{-1}L\right) = O. \]
The last equality follows from Lemma \ref{UenLequatilties}.
Thus the only non-zero elements of $MA - BM$ are on the first
subdiagonal.
\end{proof}

\begin{remark}  The generalized finite RLT
(\ref{moregeneral}) can also be written as
\begin{equation}
\label{generalinaenb} \left\{ \begin{array}{ll} \dot
a_{n}=\Bigl(F\left(LM^{-1}\right)\Bigr)_{n,n-1}-\Bigl(F\left(M^{-1}L\right)\Bigr)_{n+1,n},&
\quad 1\le n\le N,\\[2ex] \dot
b_{n}=\Bigl(F\left(M^{-1}L\right)\Bigr)_{n+1,n}-\Bigl(F\left(LM^{-1}\right)\Bigr)_{n+1,n},&
\quad 1\le n\le N-1.
\end{array}\right.
\end{equation}
\end{remark}

\begin{remark}
\label{stellinganderevorm} If we take $F(x)=1/x$ in
(\ref{moregeneral}) we obtain the first representation of the
finite RTL (\ref{hetsysteem}). The other representation
(\ref{anderevorm}) found in the literature, is obtained by taking
$F(x)=x$. To see this, we note that from (\ref{dematricesLenU}) we
have
\[M^{-1}=\left(
\begin{array}{cccccc}
1 & 0 & \cdots & \cdots & \cdots & 0\\ b_{1} & 1 & 0 & & & \vdots
\\ b_2b_1 & b_1 & 1 & 0  &  & \vdots \\ \vdots & \ddots & \ddots &
\ddots & \ddots & \vdots\\ \vdots  & & \ddots & \ddots & \ddots &
0\\ \prod\limits_{i=1}^{N-1}b_i & \cdots & \cdots & b_{N-1}
b_{N-2} & b_{N-1} & 1
\end{array} \right), \]
so that $(LM^{-1})_{i+1,i}=b_{i}(b_{i+1}+a_{i+1})$ and
$(M^{-1}L)_{i+1,i}=b_{i}(b_{i-1}+a_i)$ for $i=1,\ldots,N-1$.  So, if
$F(x)=x$, the equations (\ref{generalinaenb}) are
\[
\left\{ \begin{array}{l} \dot
a_{n}=b_{n-1}(b_{n}+a_{n})-b_{n}(b_{n-1}+a_n)=a_{n}(b_{n-1}-b_{n}),\\[2ex]
\dot
b_{n}=b_{n}(b_{n-1}+a_{n})-b_{n}(b_{n+1}+a_{n+1})=b_{n}(a_{n}-a_{n+1}+b_{n-1}-b_{n+1}),
\end{array}\right.
\]
with $b_0=b_{N}=0$.  This corresponds to (\ref{anderevorm}).
\end{remark}

\subsection{The evolution of the spectral data}
\label{section_evolution} We now describe the time evolution of
the spectral data introduced in Section \ref{detransformatie}. The
generalized RTL (\ref{moregeneral}) can then be solved with the
aid of this spectral transform which is a bijection between
$\left(\mathbb{R}^+\right)^{2N-1}$ and its image $\{ 0 < \lambda_1
< \cdots < \lambda_N, \ w_j >0, \ \sum_{j=1}^N w_j = 1\}$.  In
particular we can start from positive initial matrix data
$a_{n}(0)$ and $b_{n}(0)$ and follow the scheme in Figure
\ref{schemaplot} to find the matrix data at time $t\in
\mathbb{R}$.

Let $L(t)$ and $M(t)$ be the solutions of the generalized RTL
(\ref{moregeneral}) with
\begin{equation} \label{AtenBt}
 A(t) = -\Bigl(F\left(M^{-1}(t)L(t)\right)\Bigr)_{-} \quad \mbox{and} \quad
 B(t)=-\Bigl(F\left(L(t)M^{-1}(t)\right)\Bigr)_{-}.
\end{equation}
Let $\mathcal{L}_1$ and $\mathcal{L}_2$ be the unique solutions of the
linear differential equations
\begin{equation}
\label{leftdeL} \dot{\mathcal{L}_1}(t) =\mathcal{L}_1(t)
B(t),\qquad \mathcal{L}_1(0)=I,
\end{equation}
and
\begin{equation}
\label{leftdeL2} \dot{\mathcal{L}_2}(t) =-A(t)\mathcal{L}_2(t)
,\qquad \mathcal{L}_2(0)=I.
\end{equation}
Then
$\mathcal{L}_1$ and $\mathcal{L}_2$ are lower triangular matrices
with ones on the diagonal and therefore invertible. Furthermore
\[\frac{d}{dt}(\mathcal{L}_1L\mathcal{L}_2)=
\mathcal{L}_1BL\mathcal{L}_2+\mathcal{L}_1(LA-BL)\mathcal{L}_2
-\mathcal{L}_1LA\mathcal{L}_2= O \] and similarly
$\frac{d}{dt}(\mathcal{L}_1M\mathcal{L}_2)= O$. This means
that
\begin{equation}
\label{terugnaar0voorUenL}
\mathcal{L}_1(t)L(t)\mathcal{L}_2(t)=L(0) \qquad \mbox{and} \qquad
\mathcal{L}_1(t)M(t)\mathcal{L}_2(t)=M(0).
\end{equation}
\begin{theorem}
The eigenvalues $\lambda_j$ of the generalized eigenvalue problem
of the pair $(L,M)$ are constants of motion, which means
\begin{equation}
\label{evolutielambda_j} \dot \lambda_{j}=0,\qquad 1\le j\le N.
\end{equation}
\end{theorem}
\begin{proof} From (\ref{terugnaar0voorUenL}) it follows that
for every $z$,
\[ \mathcal{L}_1(t) (z M(t) - L(t)) \mathcal{L}_2(t) =
    z M(0) - L(0). \]
Since $\det (\mathcal{L}_1(t))=\det (\mathcal{L}_2(t))=1$ we see
that \[ \det(z M(t)-L(t)) = \det(z M(0)-L(0)). \]
This proves the theorem.
\end{proof}

To find the evolution of the spectral data $w_j$, we need the following lemma.
For a matrix $X$ we use $X_{\geq}$ to denote its upper triangular part.
So $X = X_{\geq} + X_-$.
\begin{lemma}
\label{hulptheorema} Let $\mathcal{R}$ be the unique solution of
the linear differential equations
\begin{equation}
\label{rightdeR}
\dot{\mathcal{R}}(t)=-\Bigl( F\left(L(t) M^{-1}(t) \right) \Bigr)_{\geq} \
\mathcal{R}(t),\qquad \mathcal{R}(0)=I.
\end{equation}
Then $\mathcal{R}$ is upper triangular and we have
\begin{equation}
\label{emacht}
\mathcal{L}_1(t)\mathcal{R}(t)=\rme^{-tF\left(L(0)M^{-1}(0)\right)},
\end{equation}
which gives an LR-factorization of
$\rme^{-tF\left(L(0)M^{-1}(0)\right)}$.
\end{lemma}
\begin{proof} Combining (\ref{leftdeL}), (\ref{rightdeR}), and
(\ref{AtenBt}), we obtain
\begin{eqnarray*}
\frac{d}{dt}\Bigl(\mathcal{L}_1(t)\mathcal{R}(t)\Bigr) &= &
\mathcal{L}_1(t)B(t)
\mathcal{R}(t)-\mathcal{L}_1(t)\Bigl(F\left(L(t)M^{-1}(t)\right)\Bigr)_{\geq} \
\mathcal{R}(t)
\\ & = &
-\mathcal{L}_1(t)F\left(L(t)M^{-1}(t)\right)\mathcal{R}(t).
\end{eqnarray*}
From (\ref{terugnaar0voorUenL}) we get that
\begin{eqnarray*}
F\left(L(t)M^{-1}(t)\right) &=& F\left(\mathcal{L}_1^{-1}(t)L(0)M^{-1}(0)\mathcal{L}_1(t)\right)
\\ &=&
\mathcal{L}_1^{-1}(t)F\left(L(0)M^{-1}(0)\right)\mathcal{L}_1(t)
\end{eqnarray*}
(see \cite[Definition 6.2.4]{Horn}). Hence
\[\frac{d}{dt}\Bigl(\mathcal{L}_1(t)\mathcal{R}(t)\Bigr)=-F\left(L(0)M^{-1}(0)\right)\mathcal{L}_1(t)\mathcal{R}(t),\]
which proves (\ref{emacht}) since $\mathcal{L}_1(0)\mathcal{R}(0)=I$.
\end{proof}

The theorem below gives the evolution of the spectral data $w_j$.
Moser obtained a similar expression for his spectral data in
the case of the nonrelativistic Toda lattice \cite{Moser}.
\begin{theorem}
\label{evolutiongeneral} For the generalized RTL {\rm (\ref{moregeneral})}
with positive initial data $a_n(0) > 0$, $b_n(0)>0$,  the
spectral data $w_j$ have the time evolution
\begin{equation}
\label{evogeneral}
w_{j}(t)=\frac{w_{j}(0)\rme^{-tF\left(\lambda_{j}\right)}}{\sum\limits_{k=1}^Nw_{k}(0)
\rme^{-tF\left(\lambda_{k}\right)}},\qquad 1\le j\le N.
\end{equation}
\end{theorem}
\begin{proof} Set $D=\diag(\lambda_1,\ldots,\lambda_N)$,
which is time independent. The matrix $V(t)$ contains the right
eigenvectors (with the first components equal to 1) of the
generalized eigenvalue problem of the pair $(L(t),M(t))$.  For
$t=0$ this means that $L(0)V(0)=M(0)V(0)D$. Using
(\ref{terugnaar0voorUenL}) we then obtain
\[ L(t)\mathcal{L}_2(t)V(0)=M(t)\mathcal{L}_2(t)V(0)D. \]
So the columns of the matrix $\mathcal{L}_2(t)V(0)$, which are
$\mathcal{L}_2(t)\vec{v}_j(0)$, $1\le j\le N$, are right
eigenvectors of the generalized eigenvalue problem of the pair
$(L(t),M(t))$. Because $\mathcal{L}_2(t)$ is lower triangular with
ones on the diagonal, and $\vec{v}_j(0)$ has first component equal
to $1$, we find that $\mathcal{L}_2 \vec{v}_j(0)$ has first
component equal to $1$. So we have
\begin{equation}
\label{rechtseeigenvectoren}
\vec{v}_j(t)=\mathcal{L}_2(t)\vec{v}_j(0),\qquad 1\le j\le N.
\end{equation}

Similarly we get from $U(0)^TL(0)=DU(0)^TM(0)$ and
(\ref{terugnaar0voorUenL})
\[U(0)^T\mathcal{L}_1(t)L(t)=DU(0)^T\mathcal{L}_1(t)M(t).\]
The $j$th column of the matrix $\mathcal{L}_1(t)^TU(0)$, which is
$\mathcal{L}_1(t)^T\vec{u}_j(0)$, is then a left eigenvector of
the generalized eigenvalue problem of the pair $(L(t),M(t))$ with
eigenvalue $\lambda_j$.  So it is equal to $\vec{u}_j(t)$  up to a
constant factor.  We now look at the first component of
$\mathcal{L}_1(t)^T\vec{u}_j(0)$. From Lemma \ref{hulptheorema} we
get
\[
\mathcal{L}_1(t) =
\rme^{-tF\left(L(0)M^{-1}(0)\right)}\mathcal{R}^{-1}(t).
\]
Since $\mathcal{R}(t)$ is upper triangular, we  obtain
\begin{equation} \label{firstcomp}
    \vec{u}_j(0)^{\, T} \mathcal{L}_1(t) \vec{e}_1 =
    \vec{u}_j(0)^{\, T}\rme^{-tF\left(L(0)M^{-1}(0)\right)}\left(\mathcal{R}^{-1}(t)\right)_{1,1}\vec{e}_1.
\end{equation}
From the fact that $\vec{u}_j(0)$ is a left eigenvector of $L(0) M^{-1}(0)$ with
eigenvalue $\lambda_j$, it follows that it is also a left eigenvector
of $\rme^{-tF(L(0)M^{-1}(0))}$ with eigenvalue $\rme^{-tF(\lambda_j)}$. So
\[\vec{u}_j(0)^{\, T}\rme^{-tF\left(L(0)M^{-1}(0)\right)}=\rme^{-tF\left(\lambda_j\right)}\vec{u}_j(0)^{\, T}. \]
Using this in (\ref{firstcomp}) and noting that the first
component of $\vec{u}_j(0)$ is $1$, we see that
\[\vec{u}_j(0)^{\, T}\mathcal{L}_1(t)\vec{e}_1 =
    \rme^{-tF\left(\lambda_j\right)}\left(\mathcal{R}^{-1}(t)\right)_{1,1}. \]
From this we  conclude
\begin{equation}
\label{linkseeigenvectoren} \mathcal{L}_1(t)^T\vec{u}_j(0) =
\rme^{-t F(\lambda_j)} \left(\mathcal{R}^{-1}(t)\right)_{1,1}
\vec{u}_j(t), \qquad 1\le j\le N.
\end{equation}

Combining Definition \ref{definitieW_j},
(\ref{rechtseeigenvectoren}), (\ref{linkseeigenvectoren}), and
(\ref{terugnaar0voorUenL}) we now obtain
\begin{eqnarray*}
w_j(t) & = & \left(\vec{u}_j(t)^{\,
T}M(t)\vec{v}_j(t)\right)^{-1}\\ & = &
\left(\frac{\rme^{tF\left(\lambda_j\right)}}{\left(\mathcal{R}^{-1}(t)
\right)_{1,1}}\vec{u}_j(0)^{\,
T}\mathcal{L}_1(t)M(t)\mathcal{L}_2(t)\vec{v}_j(0)\right)^{-1}\\ &
= & \left(\mathcal{R}^{-1}(t)
\right)_{1,1}\rme^{-tF\left(\lambda_j\right)}\Bigl(\vec{u}_j(0)^{\,
T}M(0)\vec{v}_j(0)\Bigr)^{-1}\\ & = & \left(\mathcal{R}^{-1}(t)
\right)_{1,1}w_j(0)\rme^{-tF\left(\lambda_j\right)}.
\end{eqnarray*}
Since the $w_j(t)$ sum to $1$ we have
\[\left(\mathcal{R}^{-1}(t)\right)_{1,1} =
    \frac{1}{\sum\limits_{k=1}^Nw_k(0)\rme^{-tF\left(\lambda_k\right)}}, \]
which completes the proof of the theorem.
\end{proof}

\begin{remark}
Standard existence and uniqueness results from the
theory of ordinary differential equations, give that the generalized
RTL (\ref{moregeneral}) with positive initial data $a_n(0)$ and $b_n(0)$
has a unique solution in a time interval containing $t=0$, so that $a_n(t)>0$,
$b_n(t)>0$.  The expressions (\ref{evolutielambda_j}) and (\ref{evogeneral})
allow us to define spectral data for all $t\in \mathbb{R}$
that remain in $\{ 0 < \lambda_1 < \cdots < \lambda_N, \ w_j >0, \
\sum_{j=1}^N w_j = 1\}$. It then follows that the system
(\ref{moregeneral}) has a unique solution valid for all time, in which the
matrix data are positive.
\end{remark}

\subsection{Example}
\label{section_example}

\begin{figure}[t]
\begin{center}
\includegraphics[scale=0.5]{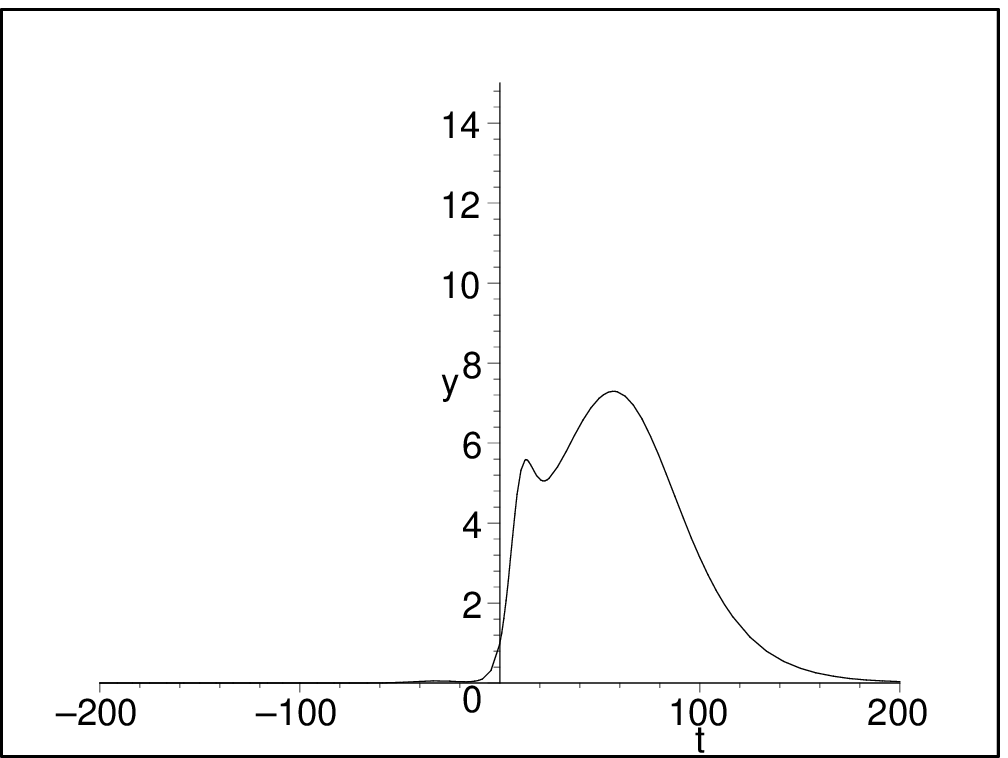}
\includegraphics[scale=0.5]{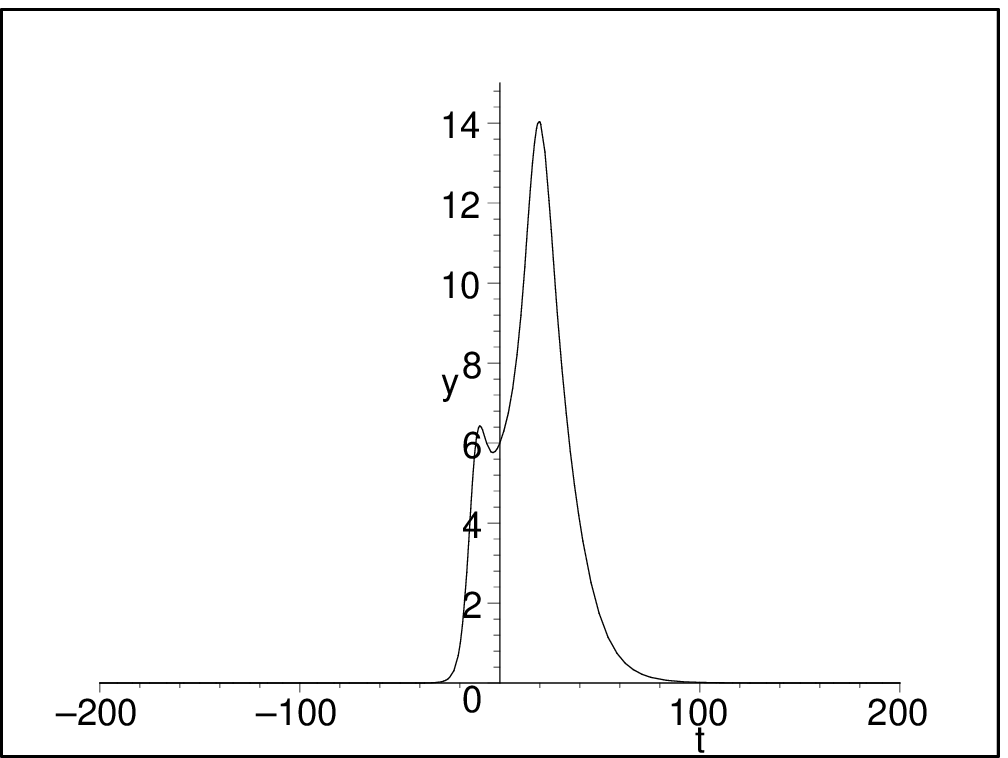}\\
\includegraphics[scale=0.5]{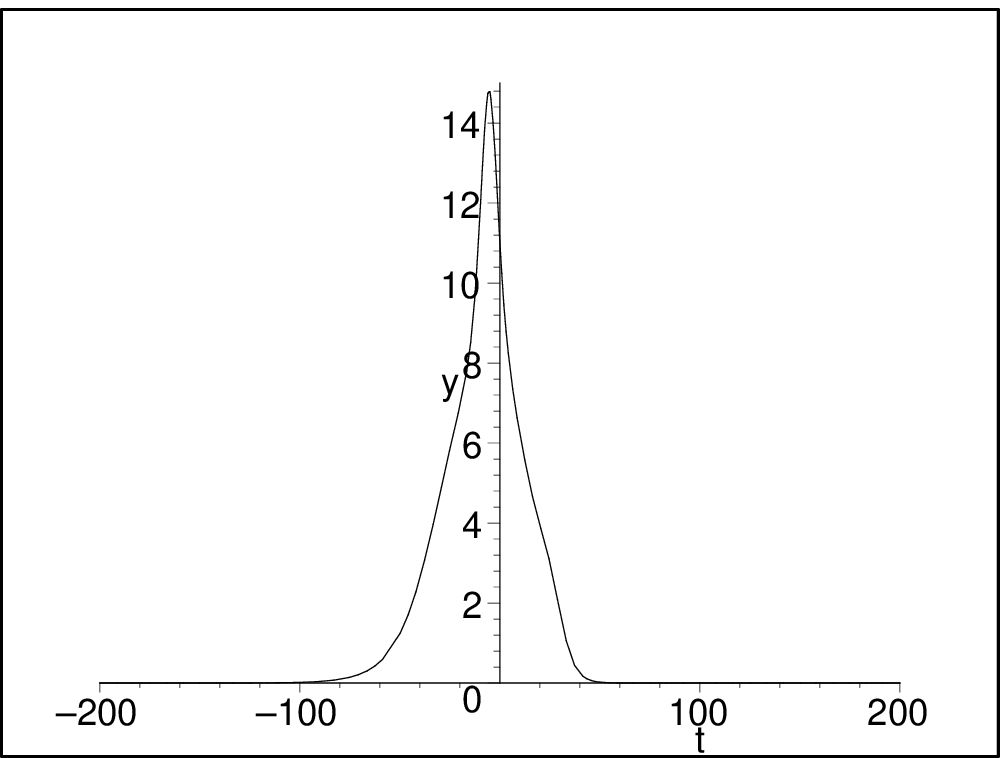}
\includegraphics[scale=0.5]{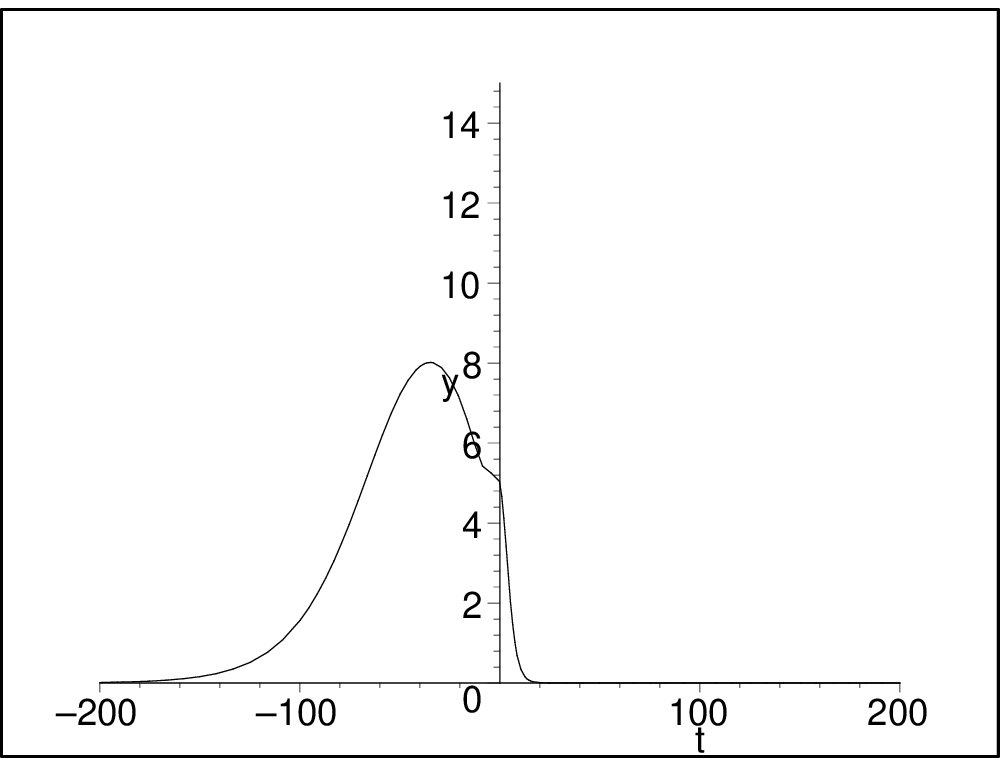}
\end{center}
\caption{\label{plotjes2}{\em Evolution of the matrix data $b_{n}$
($1\le n\le 4$).}}
\end{figure}

As an example we take $N=5$ and positive initial matrix data
\[\begin{array}{ccccccc} a_{1}(0) & = & 3, &
\qquad & b_{1}(0) & = & 1,\\ a_{2}(0) & = & 12, & \qquad &
b_{2}(0) & = & 6,\\ a_{3}(0) & = & 16, & \qquad & b_{3}(0) & = &
11,\\ a_{4}(0) & = & 7, & \qquad & b_{4}(0) & = & 5,\\ a_{5}(0) &
= & 5. & \qquad & & &
\end{array}\]
\begin{figure}[t]
\begin{center}
\includegraphics[scale=0.5]{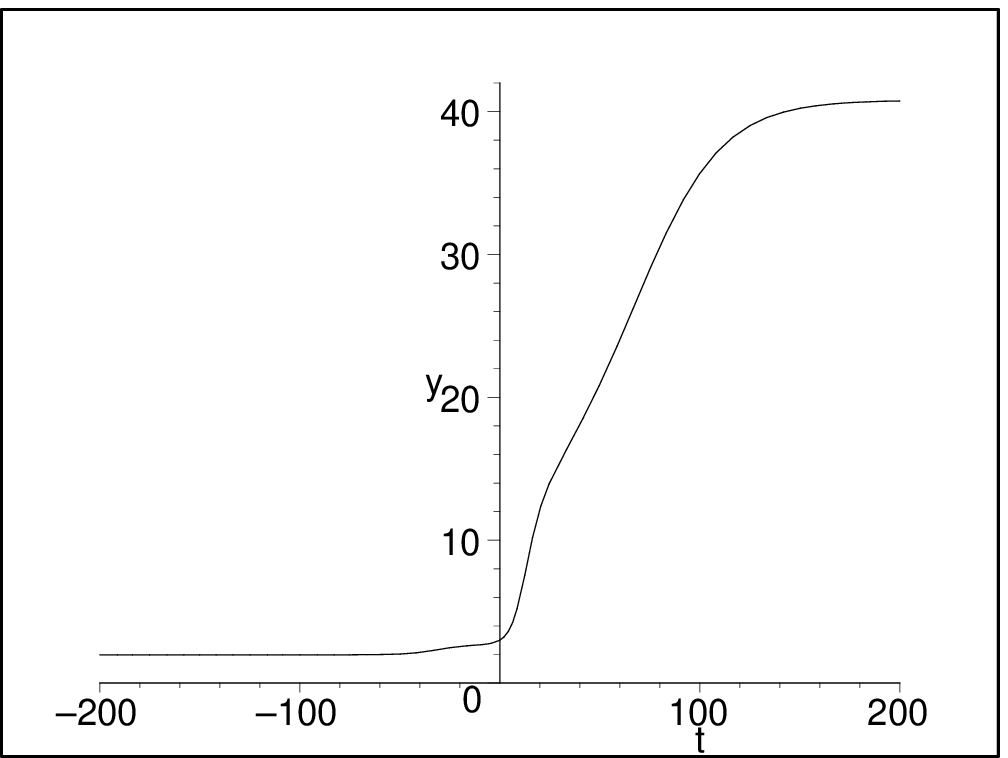}
\includegraphics[scale=0.5]{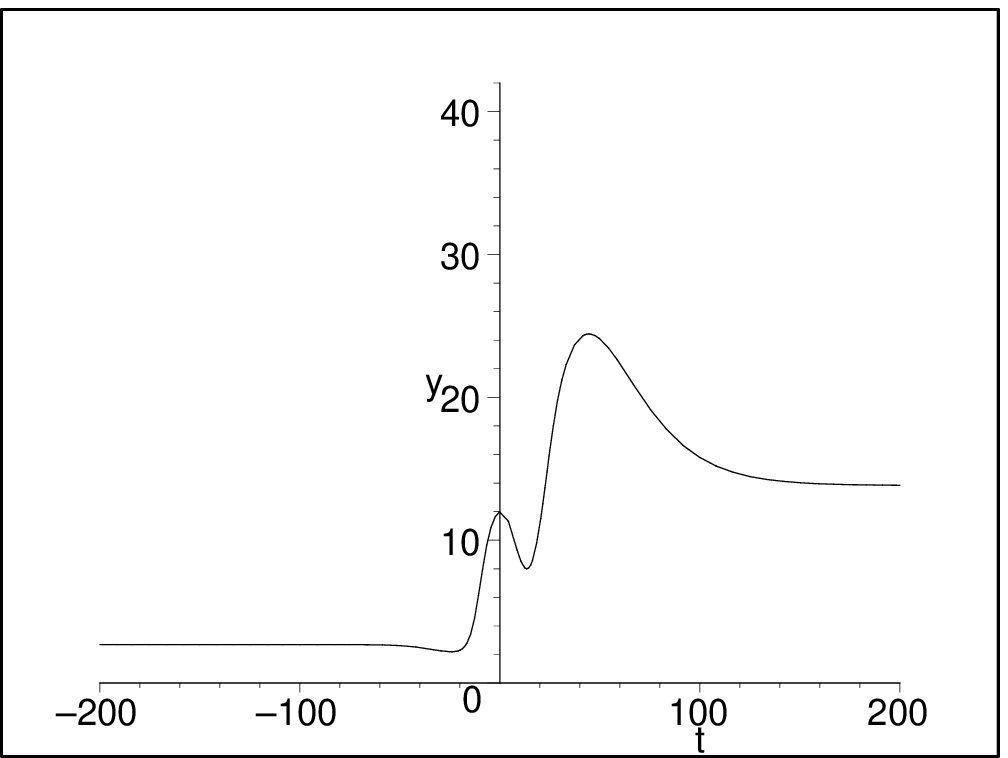}
\includegraphics[scale=0.5]{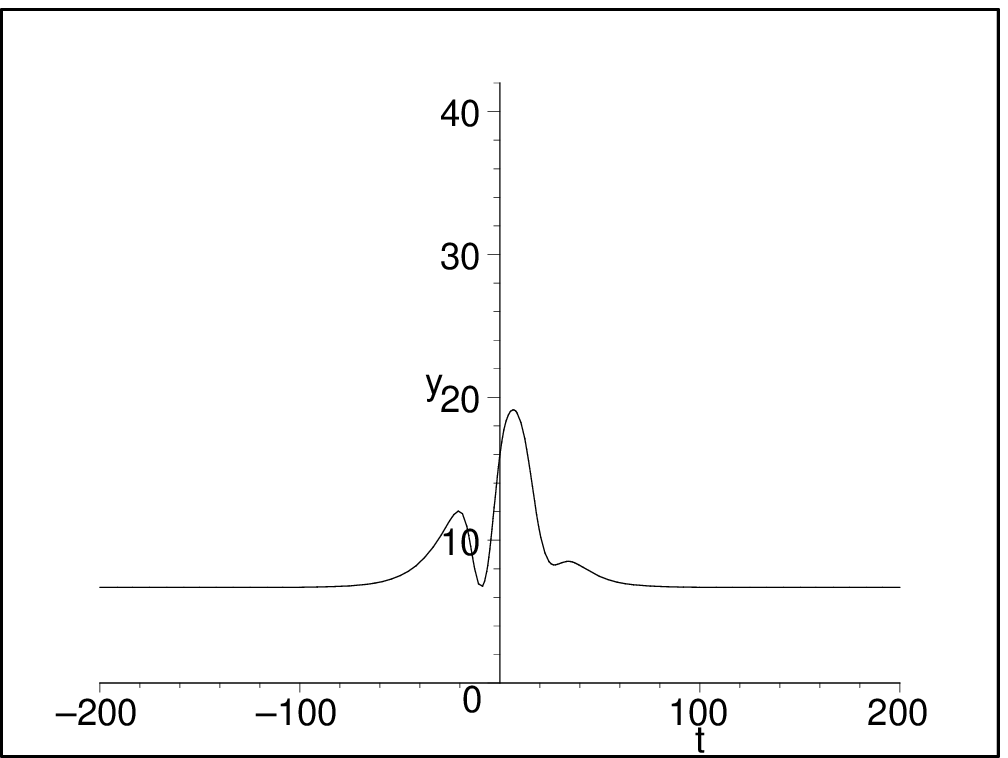}\\
\includegraphics[scale=0.5]{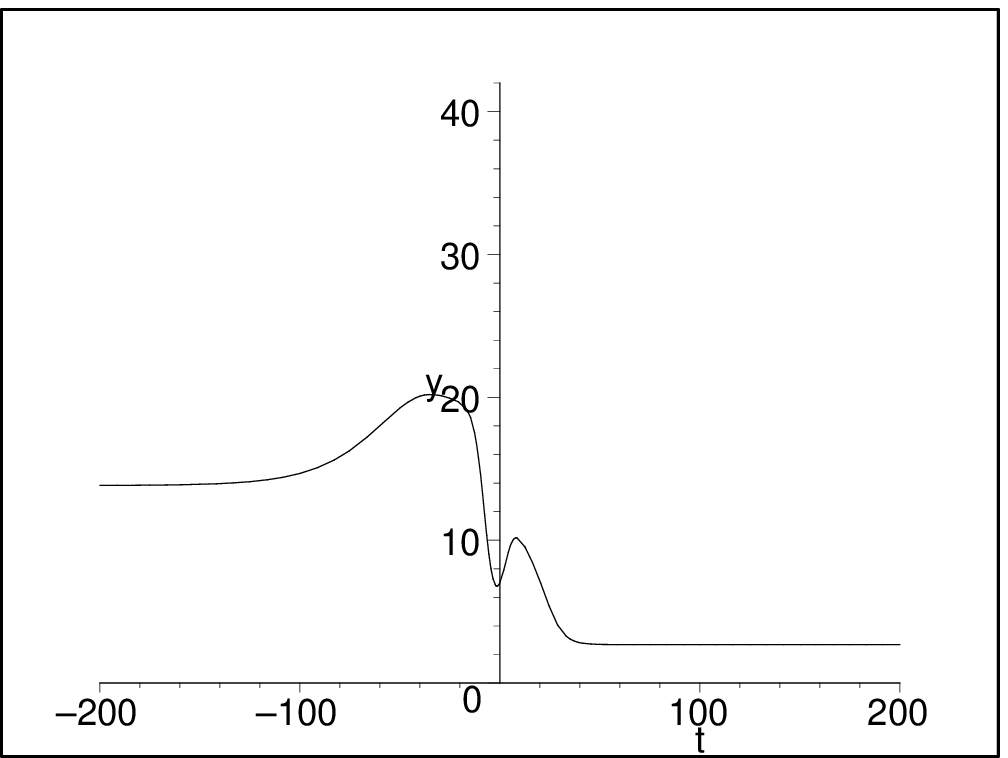}
\includegraphics[scale=0.5]{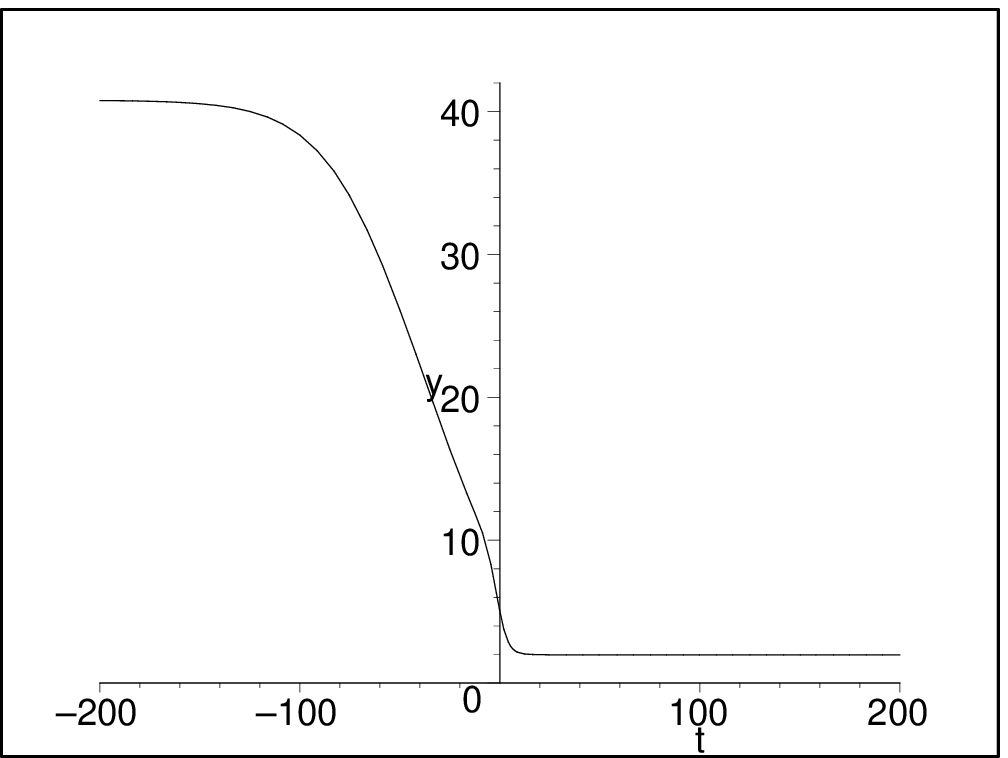}
\end{center}
\caption{\label{plotjes1}{\em Evolution of the matrix data $a_{n}$
($1\le n\le 5$).}}
\end{figure}
We obtain the corresponding spectral data, defined in Section
\ref{dsp}, by solving the generalized eigenvalue problem of the
pair $(L(0),M(0))$.  The results are
\[\begin{array}{ccccccc}
\lambda_{1} & = & 1.9812757881, & \qquad & w_1(0) & = &
0.0097186754\\ \lambda_{2} & = & 2.6941860907, & \qquad & w_2(0) &
= & 0.8409233539\\ \lambda_{3} & = & 6.6927423653, & \qquad &
w_3(0) & = & 0.0757415291\\ \lambda_{4} & = & 13.8305993379, &
\qquad & w_4(0) & = & 0.0665694128\\ \lambda_{5} & = &
40.8011964181, & \qquad & w_5(0) & = & 0.0070470286.
\end{array}\]
We consider the case that $F(x)=1/x$. Knowing
(\ref{evolutielambda_j}) and (\ref{evogeneral}) we then compute
the spectral data at several time values.  For each of these
values we construct the T-fraction (\ref{claim}) as in the proof
of Theorem \ref{existence_of_matrix_data}.  This then gives the
values of the matrix data. The results are shown in Figures
\ref{plotjes1} and \ref{plotjes2}.  We observe that the $b_n$ tend
to 0 as $t\to \pm \infty$.  It is also clear from Figure
\ref{plotjes1} that the $a_n$ have limits for $t\to \pm \infty$.
In the next section we will show that these limits are in fact
$\lambda_{6-n}$ and $\lambda_{n}$ as $t\to +\infty$, or $t\to
-\infty$, respectively.

\section{The sorting property of the generalized finite RTL}
\label{section_limits}

In this section we suppose that $F$ is strictly increasing and we
investigate the behaviour of the matrix data of the generalized
finite RTL when $t\to \pm\infty$. We will show that
\begin{equation}
\label{limietb_n}\lim_{t\to \pm \infty}b_{n}(t)=0,\qquad 1\le n\le
N-1,\end{equation} and
\begin{equation}
\label{limieta_n} \lim_{t\to +\infty}a_{n}(t)=\lambda_{n},\qquad
\lim_{t\to -\infty}a_{n}(t)=\lambda_{N-n+1},\qquad 1\le n\le N.
\end{equation}
We also establish the precise rates of convergence of
(\ref{limietb_n}) and (\ref{limieta_n}). For the SVD flow (and the
related nonrelativistic Toda lattice)  this was obtained in
\cite{Deift}.  Note that in view of the sorting property in
(\ref{limieta_n}) the choice of spectral data
$\lambda_1<\cdots<\lambda_N$ appears to be very natural. In the
limit $t \to \pm \infty$, the diagonal entries $a_n(t)$ of the
matrix $L(t)$ tend to the eigenvalues and the eigenvalues are
sorted in increasing order (as $t \to +\infty$) or in decreasing
order (as $t \to -\infty$).

The case when $F$ is strictly decreasing can be studied in the same
way as we will do for strictly increasing $F$.  In this case,
the $b_n(t)$ tend to $0$, and the $a_n(t)$ tend to the eigenvalues
in decreasing order as $t \to +\infty$, and in increasing order
as $t \to -\infty$.

Our proof of (\ref{limietb_n}) and (\ref{limieta_n}) is based on
properties of the Laurent orthogonal polynomials  $P_n$. To
indicate the time dependence we write $P_n(z, t)$. These
polynomials are defined by the recurrence relation
\begin{equation}
\label{recursiemett}
P_{n}(z,t)=(z-a_{n}(t))P_{n-1}(z,t)-b_{n-1}(t)zP_{n-2}(z,t),
\end{equation}
with  $P_0\equiv 1$ and $P_{-1}\equiv 0$. For the polynomials
$P_n(z,t)$ we find the following limits as $t\to \pm \infty$.
\begin{lemma}
\label{limietP_n}  Assume $F$ is strictly
increasing. Then we have for $1 \leq n \leq N$,
\begin{equation} \fl \qquad \qquad
\label{limietP_nuitdrukking}\lim_{t\to
+\infty}P_n(z,t)=\prod_{k=1}^n(z-\lambda_{k})\quad \mbox{and}
\quad \lim_{t\to
-\infty}P_n(z,t)=\prod_{k=N-n+1}^N(z-\lambda_{k}).\end{equation}
\end{lemma}
\begin{proof}
 First we investigate
what happens when $t\to +\infty$. In Theorem
\ref{Laurent-orthogonality} we showed that the polynomial
$P_n(z,t)$ is orthogonal with respect to the discrete measure
$\sum_{j=1}^N \frac{w_j(t)}{\lambda_j^n} \delta_{\lambda_j}$. It
is well-known that monic orthogonal polynomials minimize the
weighted $L^2$ norm, see \cite[Section 3.1]{Szego}. In our case
this means that the minimum of
\begin{equation} \label{minL2norm}
    \sum_{j=1}^N \Bigl(q_n(\lambda_j) \Bigr)^2 \frac{w_j(t)}{\lambda_j^n}
\end{equation}
taken over all monic polynomials $q_n$ of degree $n$, is attained
by $P_n(z, t)$. Take $q_n(z)=\prod_{k=1}^n(z-\lambda_{k})$, then
(since $\lambda_j >0$ and $w_j(t)>0$) we obtain for $1 \leq i \leq
N$,
\begin{eqnarray}
\nonumber \Bigl(P_n(\lambda_i,t)\Bigr)^2
\frac{w_i(t)}{\lambda_i^n} & \le & \sum_{j=1}^N
\Bigl(P_n(\lambda_j,t)\Bigr)^2 \frac{w_j(t)}{\lambda_j^n} \\ &
\leq & \nonumber \sum_{j=1}^N\Bigl(q_n(\lambda_j) \Bigr)^2
\frac{w_j(t)}{\lambda_j^n} \\ \label{ongelijkheid1bijP_n} & = &
\sum_{j=n+1}^{N}
\Bigl(\prod\limits_{k=1}^n(\lambda_{j}-\lambda_{k})\Bigr)^2
\frac{w_{j}(t)}{\lambda_j^n}.
\end{eqnarray}
From Theorem \ref{evolutiongeneral} we get for the ratio of the
spectral data $w_j(t)$ that
\begin{equation}
\label{verhoudingw_tjes} \frac{w_j(t)}{w_k(t)}=
\frac{w_j(0)}{w_{k}(0)}\rme^{-t\Bigl(F\left(\lambda_{j}\right)-F\left(\lambda_k\right)\Bigr)},\qquad
1\le j,k\le N .
\end{equation}
Since $F$ is strictly increasing, we obtain from (\ref{verhoudingw_tjes}),
\begin{equation}
\label{verhoudingw_tjeslimiet} \lim_{t\to
+\infty}\frac{w_j(t)}{w_k(t)}=0,\qquad 1\le k<j\le N.
\end{equation}
So (\ref{ongelijkheid1bijP_n}) gives for every $1 \leq i \leq N$,
\begin{equation} \fl \qquad
\label{ongelijkheidbijP_n}
\Bigl(P_n(\lambda_i,t)\Bigr)^2 \frac{w_i(t)}{\lambda_i^n}  \le
\Bigl(\prod\limits_{k=1}^n
 (\lambda_{n+1}-\lambda_{k})\Bigr)^2 \frac{w_{n+1}(t)}{\lambda_{n+1}^n}
 \Bigl(1+ {\rm o}(1)\Bigr),
 \quad \mbox{as } t\to +\infty.
\end{equation}
For $1\le i\le n$, we get from  (\ref{verhoudingw_tjeslimiet})
that $\frac{w_{n+1}(t)}{w_i(t)}$ tends to $0$ as $t \to +\infty$.
So from (\ref{ongelijkheidbijP_n}) we conclude that
\[\lim_{t\to +\infty}P_n(\lambda_{i},t)=0,\qquad 1\le i\le  n.\]
Then the first limit in (\ref{limietP_nuitdrukking}) follows,
since $P_n(z, t)$ is a monic polynomial of degree $n$ for every
$t$ . Similarly we can prove the second limit when $t\to -\infty$,
by taking
\[ q_n(z)=\prod\limits_{i=N-n+1}^N(z-\lambda_{i}) \]
in (\ref{minL2norm}).
\end{proof}

Starting from the recurrence relation (\ref{recursiemett}), the
limits of the polynomials $P_n(z, t)$ as $t \to \pm \infty$ give
rise to the limits of the matrix data (\ref{limietb_n}) and
(\ref{limieta_n}).  However, we also want to establish the precise
rate of convergence.  We then need the following lemma.
\begin{lemma}
\label{productvanbtjes} Let $0\le n\le N-1$.
Then for $1\le l\le n+1$,
\begin{equation} \fl
\label{productbtjes+} \qquad \prod_{k=1}^nb_k(t)=P_n(\lambda_l,t)
\left(\prod_{i=1,i\not=l}^{n+1}(\lambda_l-\lambda_{i})\right)
\frac{w_l(t)}{\lambda_l^n}\ \Bigl(1+{\rm o}(1)\Bigr),\quad
\mbox{as}\ t\to +\infty,
\end{equation}
and  for $N-n\le l\le N$,
\begin{equation} \fl
\label{productbtjes-} \qquad \prod_{k=1}^nb_k(t)=P_n(\lambda_l,t)
\left(\prod_{i=N-n,i\not=l}^{N}(\lambda_l-\lambda_{i})\right)
\frac{w_l(t)}{\lambda_l^n}\ \Bigl(1+{\rm o}(1)\Bigr),\quad
\mbox{as}\ t\to -\infty.
\end{equation}
\end{lemma}
\begin{proof} From Theorem \ref{Laurent-orthogonality} it
follows that for every monic polynomial $Q_n$ of degree $n$, we
have
\[\sum_{j=1}^N
P_n(\lambda_j,t) Q_n(\lambda_j) \frac{w_j(t)}{\lambda_j^n}=\prod_{k=1}^nb_k(t),
\]
where $0\le n\le N-1$. If we take $1\le l\le n+1$ and
\[Q_n(z)=\prod_{i=1,i\not=l}^{n+1}(z-\lambda_{i}),\]
then we obtain
\begin{eqnarray} \nonumber
\prod_{k=1}^nb_k(t)&=&P_n(\lambda_l,t)\left(\prod_{i=1,i\not=l}^{n+1}(\lambda_l-\lambda_{i})\right)\frac{w_l(t)}{\lambda_l^n}
\\[1ex] && \qquad \label{somvoorproductvanb_tjes} +\sum_{j=n+2}^N
P_n(\lambda_j,t)\left(\prod_{i=1,i\not=l}^{n+1}(\lambda_j-\lambda_{i})\right)\frac{w_j(t)}{\lambda_j^n}.
\end{eqnarray}
From Lemma \ref{limietP_n} we see that
\begin{equation}
\label{verschillendvan0} \lim_{t\to +\infty}P_n(\lambda_j,t) \neq
0, \qquad n+1\le j\le N.
\end{equation}
Taking $l=n+1$ in (\ref{somvoorproductvanb_tjes}) and letting $t
\to +\infty$, we get with (\ref{verschillendvan0}) and
(\ref{verhoudingw_tjeslimiet}) that
\begin{equation} \fl
\label{productbtjes+metl=n+1} \qquad
\prod_{k=1}^nb_k(t)=P_n(\lambda_{n+1},t)
\left(\prod_{i=1}^{n}(\lambda_l-\lambda_{i})\right)
\frac{w_{n+1}(t)}{\lambda_l^n}\ \Bigl(1+{\rm o}(1)\Bigr),\quad
\mbox{as}\ t\to +\infty,
\end{equation}
which is (\ref{productbtjes+}) with $l=n+1$.
So, for large $t$, $\prod\limits_{k=1}^nb_k(t)$ behaves like
$w_{n+1}(t)$ up to a constant.  Knowing this, we get
(\ref{productbtjes+}) for $1 \leq l \leq n$ from
(\ref{somvoorproductvanb_tjes}), (\ref{verschillendvan0}) and
(\ref{verhoudingw_tjeslimiet}). The proof of
(\ref{productbtjes-}) is similar.
\end{proof}

From Lemma \ref{productvanbtjes} we can easily establish the rates
of convergence of $P_n(\lambda_l,t)$, $1\le l\le n$, as $t\to
+\infty$, and of $P_n(\lambda_l,t)$, $N-n+1\le l\le N$, as $t\to
-\infty$.

Now we look at the matrix data $a_n(t)$ and $b_n(t)$. Lemmas
\ref{limietP_n} and \ref{productvanbtjes} will help us to find
their behaviour when $t\to \pm \infty$.
\begin{theorem}
\label{blimieten} For the generalized finite RTL {\rm
(\ref{moregeneral})} with $F$ strictly increasing we have for
$1\le n\le N-1$,
\begin{equation}
\label{deblimiet1} \fl \qquad b_{n}(t)  =
\frac{w_{n+1}(0)}{w_{n}(0)}\frac{\left(\lambda_{n}\right)^{n-1}}{\left(\lambda_{n+1}\right)^n}
\left(\frac{\prod\limits_{i=1}^n(\lambda_{n+1}-\lambda_i)}
{\prod\limits_{i=1}^{n-1}(\lambda_{n}-\lambda_i)} \right)^2
{\rme}^{-t\Bigl(F\left(\lambda_{n+1}\right)-F\left(\lambda_{n}\right)\Bigr)}\
\Bigl(1+{\rm o}(1)\Bigr),
\end{equation}
as $t \to +\infty$, and
\begin{eqnarray} \fl \qquad \nonumber
 b_{n}(t) & = &
\frac{w_{N-n}(0)\left(\lambda_{N-n+1}\right)^{n-1}}{w_{N-n+1}(0)\left(\lambda_{N-n}\right)^n}\left(\frac{\prod\limits_{i=N-n+1}^N(\lambda_{N-n}-\lambda_i)}
{\prod\limits_{i=N-n+2}^N(\lambda_{N-n+1}-\lambda_i)}\right)^2\\[1ex]
\fl \qquad \label{deblimiet2} & & \qquad \qquad \qquad \qquad
\qquad \qquad   \times
\rme^{t\Bigl(F\left(\lambda_{N-n+1}\right)-F\left(\lambda_{N-n}\right)\Bigr)}
\Bigl(1+{\rm o}(1)\Bigr),
\end{eqnarray}
as $t \to -\infty$.
\end{theorem}
\begin{proof} We only give the proof for the case when $t\to
+\infty$.  The proof of the case $t \to -\infty$ is similar. Let
$0\le n\le N-1$. From Lemma \ref{limietP_n} we have that
\[\lim_{t\to
+\infty}P_n(\lambda_{n+1},t)=\prod_{i=1}^{n}(\lambda_{n+1}-\lambda_{i}).\]
So taking $l=n+1$ in (\ref{productbtjes+}) we obtain
\begin{equation}
\label{tussenbijproductbtjes} \prod_{k=1}^nb_k(t)=
\left(\prod_{i=1}^{n}(\lambda_{n+1}-\lambda_{i})\right)^2
\frac{w_{n+1}(t)}{\lambda_{n+1}^n}\
\Bigl(1+{\rm o}(1)\Bigr),\quad \mbox{as } t\to +\infty.
\end{equation}
Since
\[\frac{w_{n+1}(t)}{w_{n}(t)}=\frac{w_{n+1}(0)}{w_{n}(0)}
\rme^{-t\Bigl(F\left(\lambda_{n+1}\right)-F\left(\lambda_{n}\right)\Bigr)},\]
(\ref{deblimiet1}) follows immediately from
(\ref{tussenbijproductbtjes}).
\end{proof}

This theorem implies that the matrix data $b_n(t)$ tend to $0$ as $t\to
\pm \infty$ at an exponential rate. We now use the recurrence
relation for the polynomials $P_n(z,t)$ to compute the behaviour
of $a_n(t)$ as $t\to \pm \infty$.
\begin{theorem}
\label{alimieten} Let $1\le n\le N$.  For the generalized finite
RTL {\rm (\ref{moregeneral})} with $F$ strictly increasing we have
\begin{equation} \fl \label{alimiet1}
\qquad \lambda_{n}-a_n(t)  =
\frac{\lambda_{n}}{\lambda_{n}-\lambda_{n+1}} b_{n}(t) \Bigl(
1+{\rm o}(1)\Bigr)
 -\frac{\lambda_{n}}{\lambda_{n-1}-\lambda_{n}}b_{n-1}(t)
\Bigl( 1+{\rm o}(1)\Bigr),
\end{equation}
as $t \to +\infty$, and
\begin{eqnarray} \fl \qquad \nonumber
\lambda_{N-n+1}-a_n(t) & = &
\frac{\lambda_{N-n+1}}{\lambda_{N-n+1}-\lambda_{N-n}} b_{n}(t)
\Bigl( 1+{\rm o}(1)\Bigr) \\[1ex] \qquad & & \label{alimiet2}
\qquad \qquad \quad -
 \frac{\lambda_{N-n+1}}{\lambda_{N-n+2}-\lambda_{N-n+1}}b_{n-1}(t)
\Bigl( 1+{\rm o}(1)\Bigr),
\end{eqnarray}
as $t\to -\infty$.  Recall that $b_{N}\equiv 0$, $b_0\equiv 0$.
\end{theorem}
\begin{proof} We will only prove (\ref{alimiet1}) and in the
proof we assume $2 \leq n \leq N-1$.   The recurrence relation
(\ref{recursiemett}) for the polynomials $P_n(z,t)$ gives
\begin{equation}
\label{recursiehulp}
\lambda_n-a_n(t)=\lambda_n\frac{P_n(\lambda_n,t)}{\lambda_nb_n(t)P_{n-1}(\lambda_n,t)}b_{n}(t)
+\lambda_n\frac{P_{n-2}(\lambda_n,t)}{P_{n-1}(\lambda_n,t)}b_{n-1}(t).
\end{equation}
From Lemma \ref{limietP_n} we obtain
\begin{equation}
\label{bijdea_tjesverhoudingP_n}
\frac{P_{n-2}(\lambda_n,t)}{P_{n-1}(\lambda_n,t)}=
\frac{1}{\lambda_n-\lambda_{n-1}}\ \Bigl(1+{\rm o}(1)\Bigr),
\qquad \mbox{as}\ t\to +\infty.
\end{equation}
From (\ref{productbtjes+}) we have
\[
\frac{P_n(\lambda_n,t)}{\lambda_n^n\prod\limits_{k=1}^nb_k(t)}w_n(t)=
\left(\prod_{i=1,i\not=n}^{n+1}(\lambda_n-\lambda_{i})\right)^{-1}
\ \Bigl(1+{\rm o}(1)\Bigr),\qquad \mbox{as}\ t\to +\infty.
\]
and
\[
\frac{P_{n-1}(\lambda_n,t)}{(\lambda_n)^{n-1}\prod\limits_{k=1}^{n-1}b_k(t)}w_n(t)=
\left(\prod_{i=1}^{n-1}(\lambda_n-\lambda_{i})\right)^{-1} \
\Bigl(1+{\rm o}(1)\Bigr),\qquad \mbox{as}\ t\to +\infty.
\]
Taking the ratio of the last two expressions, we obtain
\begin{equation}
\label{bijdea_tjesverhoudingP_nmetb}
\frac{P_n(\lambda_n,t)}{\lambda_nb_n(t)P_{n-1}(\lambda_n,t)}=
\frac{1}{\lambda_n-\lambda_{n+1}}\
\Bigl(1+{\rm o}(1)\Bigr),\qquad \mbox{as}\ t\to +\infty.
\end{equation}
Combining (\ref{recursiehulp}), (\ref{bijdea_tjesverhoudingP_n})
and (\ref{bijdea_tjesverhoudingP_nmetb}) gives (\ref{alimiet1}).

The cases $n=1$ and $n = N$ are simpler, because then either the
last or the first term on the right-hand side of
(\ref{recursiehulp}) equals 0.  The proof for the case $t\to
-\infty$ is similar.
\end{proof}
Since we know from Theorem \ref{blimieten} that the $b_n(t)$ tend
to $0$ as $t\to \pm\infty$, Theorem \ref{alimieten} proves the
limits (\ref{limieta_n}).  Because the $b_n(t)$ are exponentially
small as $t \to \pm\infty$, we also see that $a_n(t)$ approaches
its limit at an exponential rate.

\begin{remark}
We return to the Newtonian form (\ref{Newtoniaans}) of the finite
RTL.  In particular we look at the behaviour of the differences of
the $q_n$ as $t\to \pm\infty$.  In \cite{Brus1,Brus2} it was shown
that the change of variables
\begin{equation}\fl \qquad
\label{changeofvariables} \left\{\begin{array}{rcll} a_n & = &
\ds{\frac{h(q_{n-1}-q_n)\rme^{p_n}}{h(q_n-q_{n+1})}}, & \qquad
1\le n\le N,\\[2ex]
 b_{n} & = &
\epsilon^2\ds{\frac{\exp(q_{n}-q_{n+1}+p_{n})h(q_{n-1}-q_{n})}{h(q_{n}-q_{n+1})}},
& \qquad 1\le n\le N-1,\end{array}\right. \end{equation} where
$h(x) = \sqrt{1+\epsilon^2 \rme^x}$, leads to (\ref{anderevorm}).
This system is of the form (\ref{moregeneral}) with $F(x)=x$ (see
Remark \ref{stellinganderevorm}).  Here the initial conditions
$a_n(0)$ and $b_n(0)$ can be computed from the $q_n(0)$ and
$p_n(0)$ (or $q_n(0)$ and $\dot q_n(0)$). From
(\ref{changeofvariables}) we have
$b_n=\epsilon^2\rme^{-(q_{n+1}-q_{n})}a_n$, so that
\[ q_{n+1}(t)-q_{n}(t)=
\log \left(\frac{\epsilon^2a_{n}(\pm\infty)}{b_{n}(t)}\right)
+{\rm o}(1),\qquad \mbox{as } t\to \pm\infty. \]
 The asymptotic results (\ref{limietb_n}) and (\ref{limieta_n}) for the matrix data
then show that the differences between the $q_n$, namely
$q_{n+1}-q_{n}$, tend to $+\infty$ for $t\to \pm\infty$. From
Theorem \ref{blimieten} we also know that the $b_n$ decay at an
exponential rate as $t\to \pm \infty$.  So the $q_{n+1}-q_{n}$
have linear behaviour as $t\to \pm \infty$, where the intercept
and the slope can be established explicitly.   In particular we
get
\[ q_{n+1}(t)-q_{n}(t) = t \left(\lambda_{n+1}-\lambda_{n}\right)
+{\rm O}(1),\qquad \mbox{as } t\to +\infty,
\]
and
\[ q_{n+1}(t)-q_{n}(t) = - t\left(\lambda_{N-n+1}-\lambda_{N-n}\right)
+{\rm O}(1),\qquad \mbox{as } t\to -\infty,
\]
where the values $\lambda_1<\cdots<\lambda_N$ are obtained from
the initial conditions $a_n(0)$ and $b_n(0)$.
\end{remark}

\ack J Coussement is a Research Assistant of FWO-Vlaanderen. The
research was supported by INTAS project 00-272 and by FWO projects
G.0184.02 and G.0176.02.

\end{document}